\documentclass{amsart}

\newcommand{\R}{\ensuremath{\mathbb{R}}}

   \newcommand{\pair}[2]{\langle#1,#2\rangle}

   \newcommand{\tp}[3]{\{#1#2#3\}}
\newcommand{\tpo}[1]{\{#1\}}
     \newcommand{\tpc}[3]{\{#1,#2,#3\}}
\newcommand{\jbwst}{$JBW^*$-triple}

\newcommand{\mi}[5]{
\tp{#1}{#2}{\tp{#3}{#4}{#5}}=\tp{\tp{#1}{#2}{#3}}{#4}{#5}-\tp{#3}{\tp{#2}{#1}{#4}}{#5}+\tp{#3}{#4}{\tp{#1}{#2}{#5}}
}
\theoremstyle{remark}

\newcommand{\CC}{{\bf C}}

\newtheorem{lemma}{Lemma}[section]
\newtheorem{theorem}{Theorem}
\newcommand{\pf}{\noindent{\it Proof}.}
\newtheorem{proposition}[lemma]{Proposition}
\newtheorem{remark}[lemma]{Remark}
\newtheorem{corollary}[lemma]{Corollary}

\newcommand{\jbst}{$JB^*$-triple}

\begin{document}
\title{Contractively complemented subspaces of pre-symmetric 
spaces}
\author{Matthew Neal}
\address{Department of Mathematics, Denison University, Granville,
Ohio 43023}
   \email{nealm@denison.edu}

\author{Bernard Russo}
\address{Department of Mathematics, University of California,
Irvine, California 92697-3875}
   \email{brusso@math.uci.edu}


\date{February 4, 2008}



\begin{abstract}

In 1965, Ron Douglas proved that if $X$ is a closed subspace of an
$L^1$-space and $X$ is isometric to another $L^1$-space, then $X$ is the
range of a contractive projection on the containing $L^1$-space.
In 1977 Arazy-Friedman showed that if a subspace $X$ of $C_1$ is isometric to another $C_1$-space 
(possibly finite dimensional), then there is a contractive projection of
$C_1$    
onto $X$.
In 1993
Kirchberg proved that if a subspace $X$ of the predual of a von Neumann
algebra $M$ is isometric to the predual of another von Neumann algebra,
then there is a contractive projection of the predual of $M$ onto $X$.

We widen significantly the scope of these results by showing that
if a subspace $X$ of the predual of a $JBW^*$-triple $A$ is isometric to
the predual of another $JBW^*$-triple $B$, then there is a contractive
projection on the predual of $A$ with range $X$, as long as $B$ does not have a
direct summand which is isometric to a space of the form
$L^\infty(\Omega,H)$, where $H$ is a Hilbert space of dimension at least
two.  The result is false without this restriction on $B$.

\end{abstract}

\maketitle

\section{Introduction and background}

\subsection{Introduction}
In 1965, Douglas \cite{Douglas65} proved that the range of  a contractive
projection on an $L^1$-space is isometric to another $L^1$-space.   At the
same time, he showed the converse: if $X$ is a closed subspace of an
$L^1$-space and $X$ is isometric to another $L^1$-space, then $X$ is the
range of a contractive projection.
Both of these results were shortly thereafter extended  to $L^p$-spaces,
$1<p<\infty$ by Ando \cite{Ando66} and Bernau-Lacey \cite{BerLac74}.
The first result fails for $L^\infty$-spaces as shown by work of
Lindenstrauss-Wulbert \cite{LinWul69} in the real case  and Friedman-Russo
\cite{FriRus82} in the complex case. But not by much---the image of a contractive projection on $L^\infty$ is a $C_\sigma$-space.

Moving to the non-commutative situation, it was already known in 1978
through the  work of Arazy-Friedman \cite{AraFri78}, which gave a complete
description of the range of a contractive projection on the Schatten class
$C_1$, that the range of such a projection  is not necessarily isometric to
a space $C_1$.  However, in
   1977, Arazy-Friedman \cite{AraFri77} showed that if a subspace $X$ of
$C_p$ $1\le p\le \infty,p\ne 2$ is isometric to another $C_p$-space
(possibly finite dimensional), then there is a contractive projection of
$C_p$ onto $X$.  Moreover, in 1992, Arazy-Friedman \cite{AraFri92} also
gave a precise description of the range of a contractive projection on
$C_p$, $1<p<\infty,p\ne 2$.

Generalizing the 1978 work of Arazy-Friedman on $C_1$ to an arbitrary
noncommutative $L^1$-space, namely the predual of a von Neumann algebra,
Friedman-Russo \cite{FriRus85} showed in 1985 that the range of a
contractive projection on such a predual is isometric to the predual of a
$JW^*$-triple, that is, a weak$^*$-closed subspace of
$B(H,K)$ closed under the triple product $xy^*z+zy^*x$.
Important examples of $JW^*$-triples besides von Neumann algebras and
Hilbert spaces ($H=B(H,\CC)$) are the subspaces of $B(H)$
of symmetric (or anti-symmetric) operators with respect to an involution,
and spin factors.
Actually, the Friedman-Russo result was valid for projections acting on
the predual of a $JW^*$-triple, not just on the predual of a von Neumann algebra.

A far reaching generalization of both the 1977 work of Arazy-Friedman (in
the case $p=1$) and the 1965 work of Douglas  was given by Kirchberg
\cite{Kirchberg93}  in 1993
in connection with his work on extension properties of $C^*$-algebras.
Kirchberg proved that if a subspace $X$ of the predual of a von Neumann
algebra $M$ is isometric to the predual of another von Neumann algebra,
then there is a contractive projection of the predual of $M$ onto $X$.

In view of the result of Friedman-Russo mentioned above, it is natural to
ask if the result of Kirchberg could be extended to  preduals of $JBW^*$-triples (the axiomatic version of $JW^*$-triples),
that is,
if a subspace $X$ of the predual of a $JBW^*$-triple $M$ is isometric to
the predual of another $JBW^*$-triple $N$, then is there a contractive
projection of the predual of $M$ onto $X$?
We  show that the answer is yes as long as the predual of $N$ does not have a direct
summand which is isometric to $L^1(\Omega,H)$
where  $H$ is a Hilbert space of dimension at least two.  To see
that this restriction is necessary, one has only to consider a subspace of
 $L^1$ spanned by two or more independent standard normal random
variables. Such a space is isometric to $L^2$ but cannot be the range of a
contractive projection on $L^1$ since by the result of Douglas it would
also be isometric to an $L^1$-space, and therefore one dimensional
(consider the extreme points of its unit ball).

\subsection{Projective rigidity. The main result}

A well-known and useful result in the structure theory of operator triple
systems is the ``contractive projection principle,'' that is, the fact
that the range of a contractive projection on a \jbst\  is linearly
isometric in a natural way to another \jbst\ (Kaup, Friedman-Russo).
Thus, the category of $JB^*$-triples and contractions is stable under
contractive projections.

To put this result, and this paper,  in proper prospective, let $\mathcal
B$ be the category of Banach spaces and contractions.  We shall say that a
sub-category $\mathcal S$ of $\mathcal B$ is {\bf projectively stable} if
it has the property that whenever $A$ is an object of $\mathcal S$ and $X$
is the range of a morphism of $\mathcal S$ on $A$ which is a  projection,
then $X$ is isometric (that is,  isomorphic in $\mathcal S$) to an object in
$\mathcal S$. Examples of projectively stable categories (some mentioned already) are, in chronological
order,

\begin{enumerate}
\item $L_1$, contractions (Grothendieck 1955 \cite{Grothendieck55})
\item $L^p$, $1\le p<\infty$, contractions (Douglas 1965 \cite{Douglas65},
Ando 1966 \cite{Ando66}, Bernau-Lacey 1974 \cite{BerLac74}, Tzafriri 1969
\cite{Tzafriri69}))
\item $C^*$-algebras, completely positive unital maps (Choi-Effros 1977
\cite{ChoEff77})
\item $\ell_p$, $1\le p<\infty$, contractions (Lindenstrauss-Tzafriri 1978
\cite{LinTza78})
\item $JC^*$-algebras, positive unital maps (Effros-Stormer 1979
\cite{EffSto79})
\item $TROs$ (ternary rings of operators), complete contractions (Youngson 1983 \cite{Youngson83})
\item $JB^*$-triples, contractions (Kaup 1984 \cite{Kaup84},
Friedman-Russo 1985 \cite{FriRus85})
\item $\ell^p$-direct sums of $L^p(\Omega,H)$, $1\le p<\infty$, $H$ Hilbert space, contractions
(Raynaud 2004) \cite{Raynaud04}
\end{enumerate}

Though $C_p$ $1\le p\le \infty$ is not projectively stable, the two works of
Arazy-Friedman \cite{AraFri78} and \cite{AraFri92} deserve to be on this
list.  For a survey of results about contractive projections and their ranges in K\"othe function spaces and Banach sequence spaces, see \cite{Beata01}.

It follows immediately that  if $\mathcal S$ is projectively stable, then
so is the category $\mathcal S_*$  of spaces whose dual spaces belong to $\mathcal S$.  It should be noted that  $TROs,\ C^*$-algebras and $JC^*$-algebras are not
stable under contractive projections and $JB^*$-triples are not stable
under bounded projections.

By considering the converse of the above property, one is lead to the
following definition which is the focus of the present paper.  A
sub-category $\mathcal S$ of $\mathcal B$ is {\bf projectively rigid} if
it has the property that whenever $A$ is an object of $\mathcal S$ and $X$
is a subspace of $A$ which is isometric to an object in $\mathcal S$, then
$X$ is the range of a morphism of $\mathcal S$ on $A$ which is a
projection. Examples of projectively rigid categories (the last two
inspired this paper), are, in chronological order, 

\begin{enumerate}
\item $\ell_p$, $1<p<\infty$, contractions (Pelczynski 1960
\cite{Pelczynski60})
\item $L^p$, $1\le p<\infty$, contractions  (Douglas 1965
\cite{Douglas65}, Ando 1966 \cite{Ando66}, Bernau-Lacey 1974
\cite{BerLac74})
\item $C_p, 1\le p<\infty$, contractions (Arazy-Friedman 1977
\cite{AraFri77})
\item Preduals of von Neumann algebras, contractions (Kirchberg 1993
\cite{Kirchberg93})
\item Preduals of $TROs$, complete contractions (Ng-Ozawa 2002
\cite{NgOza02})
\end{enumerate}

The last result, by Ng and Ozawa, fails in the category of operator  
spaces with complete contractions.
Referring to Kirchberg's paper, Ng and Ozawa conjectured that  ``a
similar statement holds for $JC^*$-triples.''   While we found that this
is not true in general, we have been able to prove the following, which
in view of the counterexample mentioned earlier, is the best possible.

\begin{theorem}\label{thm:1}
Let $X$ be a subspace of the predual $A_*$ of a \jbwst\ $A$.  If $X$ is
isometric to the predual of a \jbwst, then there is a contractive
projection $P$ on $A_*$  such that $X=P(A_*)\oplus^{\ell^1}Z$, where
$Z$ is isometric to a direct sum of spaces of the form $L^1(\Omega,H)$
where  $H$ is a Hilbert space
of dimension at least two,
$P(A_*)$ is isometric to the predual of some $JBW^*$-triple with no such
$L^1(\Omega,H)$-summand, and $P(Z)=0$.
\end{theorem}

In particular, the category of preduals of $JBW^*$-triples with no
summands of the above type is projectively rigid.

As has been made clear, JB*-triples are the most natural category for  
the study of contractive projections.
It is important to note that JB*-triples are also justified as a natural
generalization of operator algebras as well as because of their connections
with complex geometry. Indeed, Kaup showed in \cite{Kaup83} that
JB*-triples are exactly those Banach spaces whose open unit ball is a  
bounded symmetric domain.  Kaup's holomorphic  
characterization of JB*-triples directly
led to the proof of the projective stability of JB*-triples in  
\cite{Kaup84} mentioned above.
Many authors since have studied the interplay between JB*-triples and  
infinite dimensional holomorphy (see \cite{Far1},\cite{U1},\cite{U2}  
for surveys).

Preduals of JBW*-triples have been called pre-symmetric spaces (\cite{Edwards06}), which explains the title of this paper,
and have been proposed as mathematical models of physical systems 
(
\cite{Friedman05}
). In this model the operations on the physical system are represented by contractive projections on the pre-symmetric space.

The authors wish to acknowledge a fruitful discussion with Timur Oikhberg.

\section{Preliminaries}

\subsection{$JBW^*$-triples}

A {\it Jordan triple system} is a complex vector space $V$ with a
{\em  triple product} $\{\cdot,\cdot,\cdot\} : V \times V \times V
\longrightarrow V$ which is symmetric and linear in the outer
variables, conjugate linear in the middle variable and satisfies
the Jordan triple identity (also called the main identity),
\[\{a,b,\{x,y,z\}\} = \{\{a,b,x\},y,z\} - \{x,\{b,a,y\},z\} +
\{x,y,\{a,b,z\}\}. \] A complex Banach space $A$ is called a
 $JB^*\mbox{\it -triple}$ if it is a Jordan triple system such that
for each $z\in A,$ the linear map $$ D(z): v\in A\mapsto
\{z,z,v\}\in A $$ is Hermitian, that is, $\|e^{it D(z)}\| = 1$ for
all $t \in \R$, with non-negative spectrum in the Banach algebra
of operators generated by $D(z)$, and $\Vert D(z)\Vert =\Vert
z\Vert ^2.$ A summary of the basic facts about JB*-triples can be
found in \cite{Russo94} and some of the references therein, such
as \cite{Kaup83},\cite{FriRus85bis}, and \cite{FriRus86bis}. The operators
$D(x,y)$ and $Q(x,y)$ are defined by $D(x,y)z=\tp{x}{y}{z}$ and
$Q(x,y)z=\tp{x}{z}{y}$, so that $D(x,x)=D(x)$ and we define $Q(x)$ to be
$Q(x,x)$.

A $JB^*\mbox{-triple}\; A$ is called a  $JBW^*\mbox{\it -triple}$
if it is a dual Banach space, in which case its predual,   denoted by $A_*$,  is unique (see  
\cite{BT} and \cite{HornMS}),
 and the triple product is separately weak*
continuous.  Elements of the predual are referred to as
normal functionals.
It follows from the uniqueness of preduals that an  
isomorphism from a JBW*-triple
onto another JBW*-triple is automatically normal, that is, w*-continuous. We will use this  
fact repeatedly in the paper.
The second dual $A^{**}$ of a $JB^*\mbox{-triple}$ is
a $JBW^*\mbox{-triple.}$   

   The $JB^*\mbox{-triples}$ form a large class of Banach spaces which
    include $C^*\mbox{-algebras,}$ Hilbert spaces, spaces of
rectangular matrices, and JB*-algebras. The triple product in a
C*-algebra $\mathcal A$ is given by
$$ \{x,y,z\} = \,\frac 12\; (xy^*z+ zy^*x). $$
In a JB*-algebra with product $x\circ y$, the triple product making it into a $JB^*$-triple  is
given by $\{x,y,z\}=(x\circ y^*)\circ z+z\circ (y^*\circ
x)-(x\circ z)\circ y^*$.  

An element $e$ in a JB*-triple $A$ is
called a {\it tripotent} if $\{e, e, e \}=e$ in which case the map
$D(e) :  A \longrightarrow A$ has eigenvalues $0,\, {1\over 2}$
and $ 1$, and we have the following decomposition in terms of
eigenspaces $$ A=A_2(e)\oplus A_1(e)\oplus A_0(e), $$ which is
called the {\it Peirce decomposition} of $A$. The ${k\over
2}$-eigenspace $A_k(e)$ is called the {\it Peirce k-space}. The
{\sl Peirce projections} from $A$ onto the Peirce k-spaces are
given by $$ P_2(e) = Q^2(e), \quad P_1(e)= 2(D(e)- Q^2(e)), \quad
P_0(e)= I -2D(e)+Q^2(e) $$ where, as noted above,  $Q(e)z = \{ e, z, e\}$ for $z\in
A$. The Peirce projections are contractive.  

Tripotents $u$ and $v$ are {\it compatible} if $\{P_k(u),P_j(v):k,j=0,1,2\}$ is a commuting family.   This holds for example if $u\in A_k(v)$ for some $k$.
For any tripotent $v$, the space $A_{2}(v)$ is a JB*-algebra under
the product $x \cdot y =\{x \,\ v \,\ y \}$ and involution
$x^{\sharp}=\{v \,\ x \,\ v \}$.
Tripotents $u,v$ are {\it orthogonal} if $u\in A_0(v)$. More generally, arbitrary elements $x,y$ are
orthogonal if $D(x,y)=0$, and we write $x\perp y$ if this is the case.
 Tripotents $u,v$ are  {\it  
collinear} if $u\in
A_1(v)$ and $v\in A_1(u)$, notation $v\top u$,  and {\it rigidly collinear} if $A_2(u)\subset
A_1(v)$ and $A_2(v)\subset A_1(u)$. 

A powerful computational tool connected with Peirce decompositons is the
so-called {\it Peirce calculus}, which states that
$$\tpc{A_k(u)}{A_j(u)}{A_i(u)}\subset A_{k-j+i}(u), $$
$$\tpc{A_0(u)}{A_2(u)}{A}=\tpc{A_2(u)}{A_0(u)}{A}=0,$$
  where it is
understood that $A_j(u)=0$ if $j\not\in\{0,1,2\}$.

In the case of a tripotent $u$ in a $JBW^*$-triple $A$ with predual $A_*$,
there is a corresponding Peirce decomposition of the normal functionals:
$A_*=A_2(u)_*\oplus A_1(u)_*\oplus A_0(u)_*$ in which $A_2(u)_*$ is
linearly spanned by the normal states of the $JBW^*$-algebra $A_2(u)$.
The norm exposed face $\{f\in A_*: f(u)=1=\|f\|\}$ is automatically a
subset of $A_2(u)_*$ and coincides with the set of normal states of
$A_2(u)_*$.

Given a JBW*-triple $A$ and $f$ in the
predual $A_*$, there is a unique tripotent $v_f \in A$, called the
{\it support tripotent} of $f$, such that $f \circ P_2(v_f) = f$
and  the restriction $f|_{A_2(v_f)}$ is a {\it faithful positive}
normal functional on the $JBW^*$-algebra $A_2(v_f)$.  It is known that for any tripotent $u$, if $f\in A_j(u)_*$ ($j=0,1,2$),  then $v_f\in A_j(u)$.   The converse is true for $j=0$ or 2 but fails in general for $j=1$ (however, see the proof of Lemma~\ref{lem:6.2}).

The set of tripotents in a $JBW^*$-triple, with a largest element
adjoined,  forms a complete lattice under the order $u\le v$ if $v-u$ is a tripotent orthogonal to $u$. This lattice is  isomorphic to various collections of
faces in the $JBW^*$-triple and its predual (\cite{EdwRut88}).
A maximal element of this lattice other than the artificial largest
element is simply called a {\it maximal tripotent}, and is the same as an
extreme point of the unit ball of the $JBW^*$-triple.  Equivalently, a
maximal tripotent is one for which the Peirce 0-space vanishes, and it is also referred to as a {\it complete tripotent}.

We shall occasionally use the joint Peirce decomposition for two
orthogonal tripotents $u$ and $v$, which states that
\[
A_2(u+v)=A_2(u)\oplus A_2(v)\oplus[A_1(u)\cap A_1(v)],
\]
\[
A_1(u+v)=[A_1(u)\cap A_0(v)]\oplus [A_1(v)\cap A_0(u)],
\]
$$A_0(u+v)=A_0(u)\cap A_0(v).$$

Let $A$ be a JB*-triple. For any $a \in A$, there is a triple functional
calculus, that is,
a triple isomorphism of the closed subtriple
   $C(a)$ generated by $a$ onto the commutative
C*-algebra $C_{0}(  \mbox{Sp}\, D(a,a)\cup\{0\} )$ of continuous functions
vanishing at zero, with the triple product $f\overline{g}h$.
Any
JBW*-triple has the propertly that it is the
norm closure of the linear span of its tripotents.   This is a consequence
of the spectral theorem in $JBW^*$-triples, which states that every
element has a representation of the form $x=\int \lambda du_\lambda$
analogous to the usual spectral theorem for self-adjoint operators, in
which $\{u_\lambda\}$ is a family of tripotents \cite[Lemma
3.1]{EdwRut88}.

For any element $a$ in a $JBW^*$-triple, there is a least tripotent,
denoted by $r(a)$ and referred to as the support of $a$, such that $a$ is
a positive element in the $JBW^*$-algebra $A_2(r(a))$ (\cite[Section
3]{EdwRut88}).

A closed subspace $J$ of a $JBW^*$-triple $A$  is an {\it ideal} if $\tp{A}{A}{J}\cup\tp{A}{J}{A}\subset J$ and a weak$^*$-closed ideal $J$ is complemented in the sense that $J^\perp:=\{x\in A: D(x,J)=0\}$ is also a weak$^*$-closed ideal and 
 $A=J\oplus J^\perp$.  A tripotent $u$ is said to be a {\it central tripotent} if $A_2(u)\oplus A_1(u)$
is a weak$^*$-closed ideal, and is hence orthogonal to $A_0(u)$
(\cite{HornMS}).  The structure theory of $JBW^*$-triples has been
well developed, using this and other concepts in \cite{HornMZ} and
\cite{HorNeh88}.

\medskip

The following lemma, \cite[Lemma 1.6]{FriRus85bis}, will be used repeatedly.
\begin{lemma}\label{lem:FR85}
If $u$ is a tripotent in a $JBW^*$-triple and $x$ is a norm one element
with $P_2(u)x=u$, then $P_1(u)x=0$.  Put another way,  $x=u+q$ where
$q\perp u$.
\end{lemma}
\subsection{Some general lemmas}

\begin{lemma}\label{lem:3.13}
Let $u_\lambda$ be a family of tripotents in a $JBW^*$-triple $B$ and
suppose $\sup_\lambda u_\lambda$ exists.
\begin{description}
\item[(a)] If $u_\lambda\perp y$ for some element $y\in B$, then
$\sup_\lambda u_\lambda\perp y$.
\item[(b)] If $u_\lambda\in B_1(t)$ for some tripotent $t$, then
$\sup_\lambda u_\lambda\in B_1(t)$.
\end{description}
   \end{lemma}
\pf\
\begin{description}
\item[(a)]

If $y\perp  u_\lambda$ for all $\lambda$, then $r(y)\perp u_\lambda$.  If
we let $z=\sup u_\lambda$ and
$
z=z_2+z_1+z_0$ be the Peirce decomposition with respect to $r(y)$,
then  by Peirce calculus,
$u_\lambda=\tp{u_\lambda}{z}{u_\lambda}=\tp{u_\lambda}{z_0}{u_\lambda}$ so
that by  Lemma~\ref{lem:FR85}, $z_0=u_\lambda+b_\lambda$ with $b_\lambda
\perp u_\lambda$. Therefore $r(z_0)\ge u_\lambda$, which implies $z\le
r(z_0)\in B_0(r(y))$ and so $z\in B_0(r(y))$ and therefore $z\perp y$.

\item[(b)]
Write $\sup u_\lambda=x_2+x_1+x_0$ with respect to $t$.   Since
$D(u_\lambda,u_\lambda)(\sup u_\lambda)=u_\lambda$, by Peirce calculus we
have $D(u_\lambda,u_\lambda)x_1=u_\lambda$ and
$D(u_\lambda,u_\lambda)x_2=D(u_\lambda,u_\lambda)x_0=0$. By (a), $x_2\perp
\sup u_\lambda$ and $x_0\perp \sup u_\lambda$ so that
$0=D(x_2,x_2)(x_2+x_1+x_0)=\tp{x_2}{x_2}{x_2}+\tp{x_2}{x_2}{x_1}$. By
Peirce calculus,
$\tp{x_2}{x_2}{x_2}=\tp{x_2}{x_2}{x_1}=0$, so that $x_2=0$.

Similarly,
$0=D(x_0,x_0)(x_2+x_1+x_0)=\tp{x_0}{x_0}{x_0}+\tp{x_0}{x_0}{x_1}$,
$\tp{x_0}{x_0}{x_0}=\tp{x_0}{x_0}{x_1}=0$, so that $x_0=0$.\qed

\end{description}

\begin{lemma}\label{lem:3.14}
If $x$ and $y$ are orthogonal elements in a $JBW^*$-triple and if $z$ is
any element, then
   $$D(x,x)D(y,y)z=\tp{x}{\tp{x}{z}{y}}{y}.$$
   In other words, $D(x,x)D(y,y)=Q(x,y)^2$ for orthogonal $x,y$.
\end{lemma}
\pf\  By the main identity,
$$
\tp{z}{y}{\tp{x}{x}{y}}=\tp{\tp{z}{y}{x}}{x}{y}-\tp{x}{\tp{y}{z}{x}}{y}+\tp{x}{x}{\tp{z}{y}{y}},
$$
and the term on the left and the first term on the right are zero by
orthogonality.\qed

\begin{lemma}\label{lem:5.1}
If $w$ is a maximal tripotent, and if  $u$ and $v$ are tripotents with
$v\in B_1(u)\cap  B_2(w)$ and $u\in B_1(w)$, then $B_1(w)\cap
B_0(u)\subset B_0(v)$.
\end{lemma}
\pf\  Let $x\in B_1(w)\cap B_0(u)$.  Then
$D(x,x)v=2D(x,x)D(u,u)v=2\tp{x}{\tp{x}{v}{u}}{u}=0$ by Peirce calculus
with respect to $w$.\qed

\section{Local Jordan multipliers}

Let $\psi:B_*\rightarrow A_*$ be a linear isometry, where $A$ and $B$ are
\jbwst s. Then $\psi^*$ is a normal contraction of $A$ onto $B$ and by a
standard separation theorem, $\psi^*$ maps the closed unit ball of $A$ onto the
closed unit ball of $B$. Let $w$ be an extreme point of the closed unit
ball of $B$.  Since
$(\psi^*)^{-1}(w)\cap \hbox{ball}\, A$ is a non-empty weak$^*$-compact
convex set, it has an extreme point $v$, and in fact $v$ is an extreme
point of the closed unit ball of $A$.

\begin{lemma}\label{lem:2.1} 
With the above notation,  $\psi^*[A_1(v)]\subset B_1(w)$  and
$P_2(w)\psi^*[A_2(v)]=B_2(w)$.
\end{lemma}
\pf\
    If $f$ is a normal state of $B_2(w)$, then $\psi(f)$ has norm one and
$\psi(f)(v)=f(\psi^*(v))=f(w)=1$ so that $\psi(f)$ is a normal state of
$A_2(v)$.  Now let $x_1\in A_1(v)$ and suppose $\psi^*(x_1)=y_2+y_1$
with $0\ne y_2\in B_2(w)$ and $y_1\in B_1(w)$. There is a normal state
of $f$ of  $B_2(w)$ such that $f(y_2)\ne 0$.  Then
$\psi(f)(x_1)=f(\psi^*(x_1))=f(y_2)\ne 0$, a contradiction since
$\psi(f)$, being a state of $A_2(v)$, vanishes on $A_1(v)$.

To prove the second statement,  let  $z\in B_2(w)$. Then
$z=\psi^*(a_2+a_1)$ with $a_j\in A_j(v)$,
and by the first statement,
$z=P_2(w)z=P_2(w)\psi^*(a_2)+P_2(w)\psi^*(a_1)=P_2(w)\psi^*(a_2)$.
     \qed

\subsection{A construction of Kirchberg}

The following lemma was proved by  Kirchberg \cite[Lemma
3.6(ii)]{Kirchberg93} in the case of  von Neumann algebras. His proof,
which is valid for $JBW^*$-algebras, is repeated here for the convenience
of the reader.

\begin{lemma}\label{lem:2.2}
Let  $T$ be a normal unital contractive linear map of a $JBW^*$-algebra
$X$ onto another $JBW^*$-algebra $Y$, which maps the closed unit ball of
$X$ onto the closed unit ball of $Y$. For a projection $q\in Y$, let
$a\in X$ be of norm one such that $T(a)=1_Y-2q$. If $c$ is the
self-adjoint part of $a$, then
\begin{description}
\item[(i)] $T(c^2)=T(c)^2$
\item[(ii)] $T(x\circ c)=T(x)\circ T(c)$ for every $x\in X$.
\end{description}
\end{lemma}
\pf\  (Kirchberg \cite[Lemma 3.6(ii)]{Kirchberg93})
With $a\in X$ such that $T(a)=1_Y-2q$, let $c=(a+a^*)/2$.  Since $T$ is a
positive unital map on $X$,  $T(c)=(T(a)+T(a^*))/2=(T(a)+T(a)^*)/2=1_Y-2q$
and by Kadison's generalized Schwarz inequality (\cite{RobYou82}),
$1_Y\ge T(c^2)\ge T(c)^2=(1_Y-2q)^2=1_Y$, which proves (i).

Define a continuous $Y$-valued bilinear form  $\tilde T$ on
$X_{\hbox{s.a.}}$ by
$$\tilde T(x,z)=T(x\circ z)-T(x)\circ T(z).$$
   By Kadison's inequality again, $\tilde T(x,x)=T(x^2)-T(x)^2\ge 0$ so that
by the Schwarz inequality for positive bilinear  forms $$\|\tilde
T(x,y)\|\le \|\tilde T(x,x)\|^{1/2}\|\tilde T(y,y)\|^{1/2}.$$  Since
$\tilde T(c,c)=0$ we have $\tilde T(c,z)=0$ for all $z\in
X_{\hbox{s.a.}}$, and (ii) follows.\qed

\medskip

With the notation of Lemma~\ref{lem:2.2}, define a {\it Jordan multiplier}
(with respect to the data $(X,Y,T)$) to  to be any  element of the set
\[
M=\{x\in X: T(x\circ z)=T(x)\circ T(z)\hbox{ for all }z\in X\}.
\]

\begin{corollary}\label{cor:2.2}
Let $\psi:B_*\rightarrow A_*$ be a linear isometry, where $A$ and $B$ are
\jbwst s.  Let $w$ be an extreme point of the closed unit ball of $B$ and
let  $v$ be an extreme point of the closed unit ball of $A$ with
$\psi^*(v)=w$.
We set $V=P_2(w)\psi^*|A_2(v)$ and note that $V$  is a normal unital
contractive (hence positive) map of $A_2(v)$ onto $B_2(w)$. Then
\begin{description}
\item[(a)]
For each projection $q\in B_2(w)$, there is an element $a\in A_2(v)$ of
norm one such that $V(a)=w-2q$.
\item[(b)]
If $c$ is the self-adjoint part of the element $a$ in (a), then
\begin{description}
\item[(i)] $V(c^2)=V(c)^2$
\item[(ii)] $V(x\circ c)=V(x)\circ V(c)$ for every $x\in A_2(v)$.
\end{description}
\end{description}
\end{corollary}
   \pf\
   Part (a) follows from Lemma~\ref{lem:2.1} and part (b) follows from Lemma~\ref{lem:2.2}.   \qed

\medskip

With the notation of Corollary~\ref{cor:2.2}, define a {\it Jordan
multiplier} (with respect to the pair of extreme points $w\in B,v\in A$
with $\psi^*(v)=w$, or more precisely,  with respect to $A_2(v)$ and $V$)
to be any  element of the set
\[
M=\{x\in A_2(v): V(x\circ y)=V(x)\circ V(y)\hbox{ for all }y\in A_2(v)\},
\]
where $V=P_2(w)\psi^*|A_2(v)$.   We shall let $s$ denote the support of
$V$, that is,
$s=\inf\{p:p\hbox{ is a projection in }A_2(v), V(p)=1_B\}$.   Note that
$s$ is a multiplier by Lemma~\ref{lem:2.2}.

\medskip

The following two lemmas could easily have been stated and proved if
$A_2(v)$ and $B_2(w)$ were replaced by arbitrary $JBW^*$-algebras and $V$
was replaced by a normal unital contraction with support $s$ mapping the
closed unit ball onto the closed unit ball. This fact will be used
explicitly in the proof of Lemma~\ref{6.85}.

In the rest of  section 3, $A$ and $B$ denote $JBW^*$-triples,
$\psi:B_*\rightarrow A_*$ is a linear isometry, and $V=P_2(w)\psi^*$,
where $w$ is a maximal tripotent of $B$.

\begin{lemma}\label{lem:11141}
Let $x\in A_2(s)$ be such that $0\le x\le s$ and $V(x)$ is a projection
$q$ in $B_2(w)$.   Then $x\in M_2(s)$.
\end{lemma}
\pf\  We have $V(2x-s)=2q-w$ and by the functional calculus, $\|2x-s\|\le
1$. Then Lemma~\ref{lem:2.2}  shows that $2x-s\in A_2(s)$ is a multiplier
with respect to $(w,v)$, hence $2x-s\in M_2(s)$ and $x\in M_2(s)$. \qed

\begin{lemma}\label{lem:2.2bis}
\begin{description}
\item[(a)] $M$ is a  unital $JBW^*$-subalgebra of $A_2(v)$.
\item[(b)] $V|M$ is a normal unital Jordan $^*$-homomorphism of $M$ onto
$B_2(w)$ satisfying
$
V(\tp{x}{y}{x})=\tp{V(x)}{V(y)}{V(x)}\hbox{ for all }x\in M, y\in A_2(v).
$
\item[(c)] $V|M_2(s)$ is a normal unital Jordan $^*$-isomorphism of
$M_2(s)$ onto $B_2(w)$.
\end{description}
\end{lemma}
\pf\  $M$ is  clearly a  weak$^*$-closed self-adjoint  linear subspace of
$A_2(v)$. To prove it is a $JBW^*$-subalgebra, it suffices to show that if $c=c^*\in
M$, then $c^2\in M$, equivalently that $\tilde V(c^2,c^2)=0$, where $\tilde V(x,y)=V(x\circ y)-V(x)\circ V(y)$.

Using the Jordan algebra identity, namely $(b\circ a^2)\circ a=(b\circ a)\circ
a^2)$, and the fact that $c$ is a self-adjoint multiplier,  we have
$V(c^2)\circ V(c^2)=V(c)^2\circ V(c)^2=V(c)\circ (V(c)\circ
V(c)^2)=V(c)\circ(V(c)\circ V(c^2))=V(c)\circ (V(c\circ c^2))=V(c\circ(c\circ
c^2))=V(c^2\circ c^2)$.  Thus $\tilde V(c^2,c^2)=V(c^2\circ
c^2)-V(c^2)\circ V(c^2)=0$, proving (a).

By the definition of multiplier, $V$ is a Jordan $^*$-homomorphism of $M$
into $B_2(w)$. To show that it is onto, let $q$ be a projection in
$B_2(w)$.  By Corollary~\ref{cor:2.2} there is a self-adjoint multiplier
$c$ with $V(c)=w-2q$ and so $q=(w-V(c))/2=V((v-c)/2)$.  By the spectral
theorem in $B_2(w)$, $B_2(w)_{\hbox{s.a.}}\subset V(M)$ proving that
$B_2(w)\subset V(M)$ and hence $B_2(w)=V(M)$.
The last statement in (b) follows from the relation $\tp{x}{y}{x}=2x\circ
(x\circ y^*)-y^*\circ x^2$.

To prove (c), note that the kernel of $V|M_2(s)$ is a $JBW^*$-subalgebra
of $M_2(s)$ and is hence generated by its projections. If it contained a
non-zero projection $p$ then we would have $V(s-p)=w$, contradicting the
fact that $s$ is the support of $V$.
\qed

\subsection{The pullback map}

\begin{remark}\label{rem:3.5}
Starting with an extreme point $w\in B$, every choice of extreme point
$v\in A$ with $\psi^*(v)=w$ determines the objects $V,s,M$. This notation will prevail throughout this section. For use in the
next three lemmas, we define $\phi:B_2(w)\rightarrow M_2(s)$ to be the
inverse of the Jordan $^*$-isomorphism $V|M_2(s)$.
\end{remark}

\begin{lemma}\label{lem:4.2}
If $u=\sup_\lambda u_\lambda$ in $B$,  where each $u_\lambda $ is a
tripotent majorized by a fixed maximal tripotent $w$, then
$\phi(u)=\sup_\lambda \phi(u_\lambda)$ in $A$.
\end{lemma}
\pf\ In $M_2(s)$, $\phi(u_\lambda)\le   \sup_\lambda \phi(u_\lambda)
\le  \phi(u)\le s$ so that $u_\lambda=V(\phi(u_\lambda))\le
V(\sup_\lambda\phi(u_\lambda))\le u\le w$ and therefore $u=\sup_\lambda
u_\lambda\le V(\sup_\lambda\phi(u_\lambda))\le u$.
Thus $u=V(\sup_\lambda\phi(u_\lambda))$ and since $u$ is a projection in
$B_2(w)$ and $\sup_\lambda\phi(u_\lambda)\ge 0$,
$\sup_\lambda\phi(u_\lambda)$ is a multiplier by Lemma~\ref{lem:11141}.
Therefore
$\phi(u)=\phi(V(\sup_\lambda\phi(u_\lambda))=\sup_\lambda\phi(u_\lambda)\le\phi(u)$,
proving the lemma.\qed

\begin{lemma}\label{lem:6.1}
Let $f$ be a normal functional on $B$ and let $w$ be a maximal tripotent
in $B$ with $v_ f\le w$, giving rise to $v,M,s$ in $A$ and
$\phi:B_2(w)\rightarrow M_2(s)$.  Recall that $v_f$ denotes the support
tripotent of $f$. Then
$
v_ {\psi(f)}=\phi(v_ f)$.
\end{lemma}
\pf\
Since $B_2(v_f)\subset B_2(w)$, $f\in B_2(w)_*$.  Thus
    $$\pair{\psi(f)}{s}=\pair{\psi
(P_2(w)_*f)}{s}=\pair{f}{P_2(w)\psi^*(s)}=f(w)=\|f\|=\|\psi(f)\|,$$ so that
$v_ {\psi(f)}\le s$ and hence $v_ {\psi(f)}\in A_2(s)$.

We also have
\[
\pair{\phi(v_ f)}{\psi(f)}
=\pair{P_2(w)\psi^*(\phi(v_ f))}{f}=\pair{v_ f}{f}=\|f\|=\|\psi(f)\|,
\]
and therefore
\begin{equation}\label{eq:11071}
\phi(v_ f)\ge v_{ \psi(f)}.
\end{equation}

Let $b=P_2(w)\psi^*(v_{\psi(f)})$ so that $\|b\|\le 1$ and
$$\pair{b}{f}=\pair{\psi^*(v_
{\psi(f)})}{f}=\pair{v_{\psi(f)}}{\psi(f)}=\|\psi(f)\|=\|f\|.$$  Thus $b$
belongs to the weak$^*$-closed face in $B$ generated by $f$ (that is,
$\{x\in B:\|x\|=1,\pair{x}{f}=\|f\|\}$) and therefore by \cite[Theorem
4.6]{EdwRut88}, $b=v_ f+c$ with $c\perp v_ f$.

We then have $v_ f+c=b=P_2(w)\psi^*(v_{\psi(f)})\le P_2(w)\psi^*(\phi(v_
f))=v_ f$, so that $c=0$ and
$P_2(w)\psi^*(v_{\psi(f)})=P_2(w)\psi^*(\phi(v_ f))$.  By
Lemma~\ref{lem:11141}, $v_{\psi(f)}\in M_2(s)$ and the result follows
since $P_2(w)\psi^*$ is one to one on $M_2(s)$. \qed

\medskip

   From the previous two lemmas, we can deduce the following lemma, which will be strengthened in Lemma~\ref{6.7}.

\begin{lemma}\label{lem:4.1}
With the above notation, if $u$ is any tripotent in $B$ and $w$ is a
maximal tripotent with $u\le w$, then $\phi(u)$ depends only on $u$ and
$\psi$.
More precisely, if $w'\ge u$ is another maximal tripotent and if $v'$ is a
maximal tripotent in $A$ with $\psi^*(v')=w'$ and if $M'$ and $s'$ are the
corresponding objects such that $P_2(w')\psi^*$ is a Jordan
$^*$-isomorphism of $M'_2(s')$ onto $B_2(w')$, and
$\phi'$ denotes $(P_2(w')\psi^*|M'_2(s'))^{-1}$, then $\phi(u)=\phi'(u)$.
\end{lemma}
\pf\
By Zorn's lemma, we may write $u=\sup_\lambda v_{ f_\lambda}$ for some
family $f_\lambda$ of normal functionals on $B$.  Writing $u_\lambda$ for
$v_{ f_\lambda}$, we have
\[
\phi(u)=\phi(\sup u_\lambda)=\sup\phi(u_\lambda)
\]
and
\[
\phi'(u)=\phi'(\sup u_\lambda)=\sup\phi'(u_\lambda).
\]
By Lemmas~\ref{lem:4.2} and ~\ref{lem:6.1},  $\phi(u_\lambda)=v_
{\psi(f_\lambda)}$ and $\phi'(u_\lambda)=v_{ \psi(f_\lambda)}$.
\qed

\begin{remark}\label{rem:3.10}
We define the {\bf pullback} of a tripotent $u\in B$ to be the element
$\phi(u)$ in Lemma~\ref{lem:4.1}. By this lemma, we may unambiguously denote it
by $u_\psi$.  Thus $u_\psi$ is the unique tripotent of $A$ such that for
any maximal tripotent $w$ majorizing $u$ and any maximal tripotent $v$ of
$A$ with $\psi^*(v)=w$, giving rise to the space of multipliers $M$ and
the support $s$ of $P_2(w)\psi^*|A_2(v)$, we have $u_\psi\in M_2(s)$ and
$P_2(w)\psi^*(u_\psi)=u$.  Note that in this situation,  $s=w_\psi$.
\end{remark}

We next improve  the last assertion in Lemma~\ref{lem:2.2bis} by replacing
   $V|M_2(s)$ by $\psi^*|M_2(s)$.

\begin{lemma}\label{lem:2.3} 
$\psi^*$ agrees with $V$ on $M_2(s)$.   Hence $\psi^*|M_2(s)$ is a normal
unital Jordan $^*$-isomorphism of $M_2(s)$ onto $B_2(w)$.
\end{lemma}
\pf\
We use the notation of Lemma~\ref{lem:2.2bis}. Since $V(s)=w$, we have
$\psi^*(s)=w+x_1$ where $x_1=P_1(w)\psi^*(s)$.  Then by
Lemma~\ref{lem:FR85}, $x_1=0$, so that $\psi^*(s)=w$.

It suffices to show that $\psi^*$ maps projections of $M_2(s)$ into
$B_2(w)$.   So let $p$ be any projection in $B_2(w)$.  Since
$V(p_\psi)=p$, we have $\psi^*(p_\psi)=p+y_1$ where
$y_1=P_1(w)\psi^*(p_\psi)$.  Since $p\le w$ and $y_1\in B_1(w)$,
$P_2(p)y_1=\tp{p}{\tp{p}{y_1}{p}}{p}=0$ by Peirce calculus with respect to $w$, so that
by Lemma~\ref{lem:FR85} $y_1\perp p$.   Similarly,
$\psi^*(s-p_\psi)=w-p-y_1$ and by Lemma~\ref{lem:FR85}, $y_1\perp w-p$.
Hence $y_1\in B_0(w)=\{0\}$. \qed

\medskip

The following lemma will be improved in Lemma~\ref{lem:description} to
include the case of the Peirce 2-space. As it stands, it extends the first
statement of Lemma~\ref{lem:2.1}.

\begin{lemma}\label{lem:01081}
Let $v$ be a tripotent in $B$.   Then
\begin{description}
\item[(a)]
$\psi^*(A_1(v_\psi))\subset B_1(v)+B_0(v)$
\item[(b)] $\psi^*(A_0(v_\psi))\subset B_0(v)$.
\end{description}
\end{lemma}
\pf\
Let $f$ be a normal state of $B_2(v)$.   Then
$\pair{\psi(f)}{v_\psi}=f(v)=1=\|f\|=\|\psi(f)\|$ so that $\psi(f)$ is a
normal state of $A_2(v_\psi)$ and  hence $\psi[B_2(v)_*]\subset A_2(v_\psi)_*$.

Now if $x\in A_1(v_\psi)$ and $f\in B_2(v)_*$ is arbitrary,
$\pair{f}{\psi^*(x)}=\pair{\psi(f)}{x}=0$ and  therefore $\psi^*(x)\in
B_1(v)+B_0(v)$. This proves (a).

Now let $x\in A_0(v_\psi)$ and suppose $\|x\|=1$.   Then $\|v_\psi\pm
x\|=1$ and therefore by Lemma~\ref{lem:2.3}
\begin{eqnarray*}
\|v\pm P_2(v)\psi^*(x)\|&=& \|P_2(v)\psi^*(v_\psi)\pm P_2(v)\psi^*(x)\|\\
&\le&\|\psi^*(v_\psi)\pm \psi^*(x)\|=\|\psi^*(v_\psi\pm x)\|\le 1
\end{eqnarray*}
and since $v$ is an extreme point of the unit ball of $B_2(v)$, we have
$P_2(v)\psi^*(x)=0$.  We now have $\|v+P_1(v)\psi^*(x)+P_0(v)\psi^*(x)\|=\|v+\psi^*(x)\|=\|\psi^*(v_\psi+x)\|\le 1 $
and by Lemma~\ref{lem:FR85}, $P_1(v)\psi^*(x)=0$. \qed

\medskip

\begin{lemma}\label{6.85}
Suppose $\psi^*(x)=v$ for a tripotent $v \in B$ and an element $x\in A$
with $\|x\|=1$. Then $x=v_{\psi}+q$ for some $q \perp v_{\psi}$
\end{lemma}
\pf\   Let $w$ be a maximal tripotent of $B$ majorizing $v$ and let $v'$ be a
maximal tripotent of $A$ with $\psi^*(v')=w$.

If $z\in A_2(v_\psi)$, then
$z=\tp{v_\psi}{\tp{v_\psi}{z}{v_\psi}}{v_\psi}$. Since $v_\psi$ is a
multiplier with respect to $A_2(v')$, for all $c\in A_2(v')$ we have
$P_2(w)\psi^*(v_\psi\circ c)=v\circ P_2(w)\psi^*(c)$ .  Using this and the
general formula $\tp{z}{y}{z}=2z\circ(z\circ y^*)-y^*\circ z^2$ we obtain
$P_2(w)\psi^*\tp{v_\psi}{z}{v_\psi}=\tpc{v}{P_2(w)\psi^*(z)}{v}$. For the
same reason,
$P_2(w)\psi^*(z)=\tpc{v}{P_2(w)\psi^*\tp{v_\psi}{z}{v_\psi}}{v}=\tp{v}{\tpc{v}{P_2(w)\psi^*(z)}{v}}{v}\in
B_2(v)$, proving that $P_2(w)\psi^*[A_2(v_\psi)]\subset B_2(v)$.
In fact, $P_2(w)\psi^*[A_2(v_\psi)]=B_2(v)$, since if $p$ is any
projection in $B_2(v)$, then $p_\psi\le v_\psi$, so that $p_\psi\in
A_2(v_\psi)$ and $P_2(w)\psi^*(p_\psi)=p$.

Decomposing $x=x_{2}+x_{1}+x_{0}$ with respect to $v_{\psi}$, we notice
that by Lemma~\ref{lem:01081}, $P_2(v)\psi^*(x_2)=v$,
and  since $P_2(v)\psi^*$ is a contractive unital, hence positive, hence
self-adjoint map of $A_2(v_\psi)$ onto $B_2(v)$, $P_2(v)\psi^*(x_2')=v$
where $x_2'$ is the self-adjoint part of $x_2$ in $A_2(v_\psi)$.

Now $x_2'$ is a norm one self-adjoint element of the $JBW^*$-algebra
$A_2(v_\psi)$ which $P_2(v)\psi^*$ maps to the identity $v$ of $B_2(v)$.
Thus by Lemma~\ref{lem:2.2}, we see that $x_2'$ is a multiplier with
respect to $A_2(v_\psi)$.

We show next that  $v_\psi$ is the support of the map $P_2(v)\psi^*$.
Let $p\le v_\psi$ be a projection with $P_2(v)\psi^*(p)=v$.  Then
$P_2(w)\psi^*(p)
=v$, so that by Lemma~\ref{lem:2.2bis}, $p\in M_2(s)$ and since
$P_2(w)\psi^*$ is one-to-one there, $p=v_\psi$.

Now, since $v_\psi$ is the support of the map $P_2(v)\psi^*$, it is a
multiplier with respect to $A_2(v_\psi)$, and we have $x_2'=v_\psi$ by
Lemma~\ref{lem:2.2bis} (replacing   $B_2(w)$ there by  $B_2(v)$ and
$A_2(v)$ by $A_2(v_\psi)$).

Thus $x_2=x_2'+ix_2''=v_\psi+ix_2''$ with $x_2''$ self-adjoint and by a
familiar argument, if $x_2''\ne 0$, then $\|x_2\|=\|v_\psi+ix_2''\|>1$, a
contradiction.  We now have $x_2=v_\psi$ and the proof is completed by
applying Lemma~\ref{lem:FR85} to show that $x_1=0$. \qed

\begin{remark}\label{rem:4.17} 
Suppose $x$ lies in $B$ and let $w$ be a maximal tripotent majorizing
$r(x)$.  The Jordan $^*$-isomorphism $(\psi^*|M_2(s))^{-1}$ of $B_2(w)$ onto $M_2(s)$ carries
$B_2(r(x))$ onto $M_2((r(x)_\psi)$.  We let $x_\psi$ denote the image of
$x$ under this map so that $\psi^*(x_\psi)=x$. This is an extension of the pullback of a tripotent in Remark~\ref{rem:3.10}.
\end{remark}

The following lemma shows  that $x_\psi$ may be computed using any maximal tripotent $w$ for which $x\in B_2(w)$, that is, $r(x)$ need not be majorized by $w$.   This fact will be critical in the proofs of Theorem 2 and elsewhere in this paper (for example, Lemmas~\ref{lem:6.4} and \ref{6.8}). 

\begin{lemma}\label{6.7}
Suppose $x$ is an element in $B_{2}(w)$, where $w$ is a  maximal tripotent
not necessarily majorizing $r(x)$. Let $M$ be the space of multipliers
corresponding to a choice of maximal tripotent $v$ such that
$\psi^*(v)=w$. 
Then $x_{\psi}=(\psi^*|M_2(w_\psi))^{-1}(x)$.
\end{lemma}
\pf\
We shall consider first the case that $x=u$ is a tripotent.  Let $w'$ be a
maximal tripotent majorizing $u$, so that by Lemma~\ref{lem:01081},
$\psi^*|M'_2(s')$ is a Jordan $^*$-isomorphism onto $B_2(w')$,
$u_\psi=(\psi^*|M'_2(s'))^{-1}(u)$ and let $m$ denote $(\psi^*|M_2(s))^{-1}(u)$.
Here of course, $s=w_\psi$ and $s'=w'_\psi$.

Since $\psi^*(m)=u$, by Lemma~\ref{6.85},  $m=u_\psi+q$ with $q\perp
u_\psi$. Furthermore, $\psi^*(q)=0$.

Note that since $m$ and $u_\psi$ are tripotents,  cubing the relation
$m=u_\psi+q$ shows that $q$ is also a tripotent.
We claim that $u_\psi$ and $q$ belong to $A_2(s)$.  First of all, since
$m\in A_2(s)$, we have $A_2(m)\subset A_2(s)$ and since $u_\psi\le m$ and
$q\le m$,
   $u_\psi,q\in A_2(m)\subset A_2(s)$, proving the claim.

   It remains to show that $q=0$.   To this end, note first that in
$A_2(s)$, $\tp{q}{q}{s}=q\circ q^*$ and
   $\tp{m}{q}{s}=m\circ q^*$. Using this and the fact that $m$ is a
multiplier, with $V=P_2(w)\psi^*$, we have
   \[
V(q\circ q^*)=V\tp{q}{q}{s}=V\tp{m}{q}{s}=V(m\circ q^*)=V(m)\circ V(q^*)
   =V(m)\circ V(q)^*=0.
   \]

   Now we have $V(s-q\circ q^*)=w$ so that by Lemma~\ref{lem:11141},
$s-q\circ q^*\in M_2(s)$.  Thus $q\circ q^*\in M_2(s)$ and since $V$ is
bijective on $M_2(s)$, $q\circ q^*=0$ and $q=0$.

Having proved the lemma for tripotents, we now let $x=\int \lambda du_\lambda$ be the spectral decomposition of $x$
and let $w'$ be a maximal tripotent majorizing $r(x)$.   Then for any
spectral tripotent $u_S$, we have $u_S\in B_2(w)$ and $u_S\le w'$ so that
by the special case just proved, $(u_S)_\psi=\phi(u_S)$ where
$\phi=(\psi^*|M_2(s))^{-1}$.
Approximating $x$ by  $y=\sum \lambda_iu_{S_i}$, we have
\[
y_\psi=(\psi^*|M'_2(s'))^{-1}(\sum \lambda_i u_{S_i})=\sum
\lambda_i(\psi^*|M'_2(s'))^{-1}(u_{S_i})=\sum
\lambda_i\phi(u_{S_i})=\phi(y),
\]
which completes the proof, as the maps in question are continuous.
   \qed

\medskip

\begin{remark}\label{rem:4.19}
We will henceforth refer to elements $x_{\psi}$ as multipliers without
specifying the  Peirce 2-space containing $x$.
By embedding two orthogonal elements  $x$ and $y$ of $B$ into $B_2(w)$ for
some maximal tripotent $w$, it follows that $x_\psi\perp y_\psi$.  This fact
will be used explicitly in the rest of this paper.
\end{remark}

\section{Analysis of tripotents and pullback of the Peirce
1-space}\label{sect:3}


Our next goal is to prove, in the case where $B$ has no summand  
isometric to $L^\infty(\Omega,H)$, that if $u$ is any tripotent in $B_1(w)$ for some maximal tripotent $w$, then $u_\psi\in  
A_1(w_\psi)$.
This will be achieved in this section  (see Theorem~\ref{thm:2} below) after some analysis of tripotents
in a $JBW^*$-triple.

\subsection{Rigid collinearity}

\begin{proposition}\label{prop:3.2}
If $u$ is a tripotent in $B_1(w)$ and $w$ is a maximal tripotent, then the element
$2\tp{u}{u}{w}$, which we shall denote by $w_u$, is a tripotent in $B_2(w)$
which is collinear to $u$ and $\le w$. Moreover, $u$ and $w_u$ are rigidly
collinear.
\end{proposition}

The proof will be contained in Lemmas~\ref{lem:1} to ~\ref{lem:5.5} in
which the standing assumption is that $w$ is a maximal tripotent in $B$
and $u$ is a tripotent in $B_1(w)$. This proposition was proved in  \cite[Lemma
2.5]{HornMZ} for $w$ not  necessarily maximal but under the additional
assumption that $B_2(u)\subset B_1(w)$, which follows from the maximality
of $w$.  On the other hand, Lemmas~\ref{lem:2} and ~\ref{lem:3}
are stated here with an assumption weaker than maximality and will be used in
that form later on.
   For this reason, we include the proof of Proposition~\ref{prop:3.2} here.

\begin{lemma}\label{lem:1}
If $w$ is maximal, then 
$B_2(u)\subset B_1(w)$.
\end{lemma}
\pf\ If $x\in B_2(u)$, then $x=P_2(u)x=\tp{u}{\tp{u}{x}{u}}{u}\in B_1(w)$
by Peirce calculus with respect to $w$ and the maximality of $w$.\qed
\begin{lemma}\label{lem:2}
   If $\tp{u}{w}{u}=0$ (in particular, if $w$ is maximal), then $w_u\in
B_1(u)$.
\end{lemma}
\pf\ By the main identity,
\[
\tpo{wuu}=\tp{w}{u}{\tpo{uuu}}=\tp{\tpo{wuu}}{u}{u}-\tp{u}{\tpo{uwu}}{u}+\tp{u}{u}{\tpo{wuu}}
\]
and the middle term is zero by assumption.   Hence
$$w_u/2=\tp{w_u}{u}{u}/2+\tp{u}{u}{w_u}/2=\tp{u}{u}{w_u}.\qed
$$
\begin{lemma}\label{lem:3}
   If $\tp{u}{w}{u}=0$  (in particular, if $w$ is maximal), then $w_u$ is a nonzero tripotent and $w_u\le w$.
\end{lemma}
\pf\
Clearly $w_u$ is non-zero since $u \ne 0$ does not lie in $B_{0}(w)$.  
By the main identity,
\[
\tp{u}{u}{\tpo{www}}=\tp{\tpo{uuw}}{w}{w}-\tp{w}{\tpo{uuw}}{w}+\tp{w}{w}{\tpo{uuw}}
\]
so that
\[
\tp{w}{\tpo{uuw}}{w}=2\tp{\tpo{uuw}}{w}{w}-\tpo{uuw}=2\tpo{uuw}-\tpo{uuw}=\tpo{uuw}
\]
proving that $w_u$ is a self-adjoint element of $B_2(w)$.

It remains to show that $w_u$ is an idempotent in $B_2(w)$. To this end
use the main identity to obtain
\begin{eqnarray}\nonumber
\tpo{w_uww_u}&=&2\tp{w_u}{w}{\tpo{uuw}}\\\label{eq:12042}
&=&2\left[\tp{\tpo{w_uwu}}{u}{w}-\tp{u}{\tpo{ww_uu}}{w}+\tp{u}{u}{\tpo{w_uww}}
\right].
\end{eqnarray}
Since $w_u\in B_2(w)$, the third term in the bracket on the right is equal to
$\tp{u}{u}{w_u}=w_u/2$ by Lemma~\ref{lem:2}.  It remains to show  
that the first two terms on
the right side of (\ref{eq:12042}) cancel out. In the first place,
by the main identity
\begin{eqnarray*}
u/2&=&\tp{u}{u}{\tp{w}{w}{u}}\\
&=&\tp{\tp{u}{u}{w}}{w}{u}-\tp{w}{\tp{u}{u}{w}}{u}+\tp{w}{w}{\tp{u}{u}{u}}\\
&=&\tp{\tp{u}{u}{w}}{w}{u}-\tp{w}{\tp{u}{u}{w}}{u}+u/2,
\end{eqnarray*}
so that $\tp{\tp{u}{u}{w}}{w}{u}=\tp{w}{\tp{u}{u}{w}}{u}$, that is, 
$
\tp{w}{w_u}{u}=\tp{w_u}{w}{u}.
$

On the other hand, by the main identity,
\begin{eqnarray*}
\tp{u}{w}{w_u}&=&2\tp{u}{w}{\tp{w}{u}{u}}\\
&=&2[\tp{\tp{u}{w}{w}}{u}{u}-\tp{w}{\tp{w}{u}{u}}{u}+\tp{w}{u}{\tp{u}{w}{u}}]\\
&=&2[u/2-\tp{w}{w_u}{u}/2+0]=u-\tp{w}{w_u}{u},
\end{eqnarray*}
and it now follows that $\tp{u}{w}{w_u}=\tp{w}{w_u}{u}=u/2$, proving that
the first two terms in (\ref{eq:12042}) do cancel out.\qed

\begin{lemma}\label{lem:4}
If $w$ is maximal, then $B_2(u)\subset B_1(w_u)$.
\end{lemma}
\pf\
By the joint Peirce decomposition and Lemma~\ref{lem:1},
\[
B_2(u)\subset B_1(w)=B_1(w_u)\cap B_0(w-w_u)+B_1(w-w_u)\cap B_0(w_u).
\]
Now
\[
2D(u,u)(w-w_u)=w_u-2D(u,u)w_u=w_u-w_u=0
\]
so that $u\perp (w-w_u)$ and therefore $B_2(u)\perp (w-w_u)$.  This shows
that $B_2(u)\subset  B_1(w_u)\cap B_0(w-w_u)\subset B_1(w_u)$.\qed

\begin{lemma}\label{lem:5.5}
If $w$ is maximal, then $B_2(w_u)\subset B_1(u)$; (this completes the proof of the rigid
collinearity of $w_u$ and $u$).
\end{lemma}
\pf\   Let $x\in B_2(w_u)$.  By Lemma~\ref{lem:2} and Peirce calculus with respect to $u$,  $\tpc{w_u}{P_0(u)x}{w_u}\in B_2(u)$ and by Lemma~\ref{lem:4}, $B_2(u)\subset B_1(w_u)$. 
By compatibility of $u$ and $w_u$, $P_0(u)x\in B_2(w_u)$ and by Peirce calculus with respect to $w_u$,  $P_0(u)x=\tp{w_u}{\tpc{w_u}{P_0(u)x}{w_u}}{w_u}\in B_1(w_u)$.
Hence
$P_0(u)x\in B_1(w_u)\cap B_2(w_u)=0$.
   On the other hand, by Lemma~\ref{lem:4}, $P_2(u)x\in B_1(w_u)$ so that
$P_2(u)x=0$ also.\qed

\medskip

The next two lemmas give important properties of $w_u$.

\begin{lemma}\label{lem:5} 
If $u\in B_1(w)$ and $w$ is maximal, then
$B_1(w)\cap B_0(u)\subset B_0(w_u)$.  In particular, if   $w_u=w$, then
$u\top w$ and $u$ is maximal.
\end{lemma}
\pf\
The first statement holds by Lemma~\ref{lem:5.1}.

Suppose now that $w=w_u$ so that $u\top w$.  We shall show that
$B_0(u)\subset B_0(w)$, which implies the second assertion.   By
Lemma~\ref{lem:5.5}, $B_2(w)=B_2(w_u)\subset B_1(u)$.
If $x\in B_0(u)=[B_0(u)\cap B_2(w)]+[B_0(u)\cap B_1(w)]$, say $x=x_2+x_1$
with respect to $w$, then by the first statement, $x_1\in
B_0(w_u)=B_0(w)=0$.  On the other hand, $x_2\in B_2(w)\cap B_0(u)\subset
B_1(u)\cap B_0(u)$, so $x_2=0$.\qed

\begin{lemma}\label{lem:7}
Suppose that $u_1,u_2\in B_1(w)$ with $w$ a maximal tripotent in $B$.  If
$u_1\le u_2$ then $w_{u_1}\le w_{u_2}$ and $w_{u_2-u_1}=w_{u_2}-w_{u_1}$.
\end{lemma}
\pf\
If $u_1\le u_2$, then $u_2-u_1\perp u_1$, $\tp{w}{u_1}{u_2}=\tp{w}{u_1}{u_1}$ and
\begin{eqnarray*}
w_{u_2-u_1}&=&2\tpc{w}{u_2-u_1}{u_2-u_1}\\
&=&2\tpc{w}{u_2-u_1}{u_2}-2\tpc{w}{u_2-u_1}{u_1}\\
&=&2\tp{w}{u_2}{u_2}-2\tp{w}{u_1}{u_1}-0=w_{u_2}-w_{u_1}.
\end{eqnarray*}

On the other hand, if $v_1\perp v_2$, then
by
Lemma~\ref{lem:5},  $v_2\perp w_{v_1}$ and since $w_{v_1}\perp w-w_{v_1}$,
\begin{eqnarray*}
\tp{w_{v_1}}{w_{v_1}}{w_{v_2}}&=&2\tp{w_{v_1}}{w_{v_1}}{\tp{w}{v_2}{v_2}}\\
&=&2\tp{\tp{w_{v_1}}{w_{v_1}}{w}}{v_2}{v_2}-2\tp{w}{\tp{w_{v_1}}{w_{v_1}}{v_2}}{v_2}+2\tp{w}{v_2}{\tp{w_{v_1}}{w_{v_1}}{v_2}}\\
&=&2\tp{\tp{w_{v_1}}{w_{v_1}}{w_{v_1}}}{v_2}{v_2}-0+0=2\tp{w_{v_1}}{v_2}{v_2}=0.
\end{eqnarray*}

Combining the results of the previous two paragraphs, if $u_1\le u_2$,
then $u_1\perp u_2-u_1$,
$w_{u_2-u_1}\perp w_{u_1}$, $(w_{u_2}-w_{u_1})\perp w_{u_1}$ so that
$w_{u_1}\le w_{u_2}$.\qed

\subsection{Central tripotents}

\begin{lemma}\label{lem:horn}
Let $w$ be a maximal tripotent of $B$ and suppose that $v$ is a
tripotent $\le w$, $u$ is a tripotent in $B_1(w)$ and $u\top v$.    Then
either
$B_1(w)\cap B_1(u)\cap B_0(v)\ne 0$ or $u$ is a central tripotent in $B$.
\end{lemma}
\pf\
If $v=w$ then the result follows from Lemma \ref{lem:5} so we assume 
$v\ne w$.
Suppose first  that $B_1(w)\cap B_1(u)\cap B_0(v) =0$ and  let $e\in
B_1(v)\cap B_1(w-v) \subseteq B_{2}(w)$ be a tripotent.   We shall show
that $e=0$ from which it will
follow  that $u$ is central.

We first note that $D(u)(w-v)=D(u)w-D(u)v=0$ so $w-v \in B_{0}(u)$. By Peirce calculus,  $\tpc{u}{e}{w-v}\in  B_1(w)\cap B_1(u)\cap B_0(v) =0$ and
$\tp{u}{e}{v}\in B_2(v)\cap B_1(w)\subset B_2(w)\cap B_1(w)=0$, so that
$\tp{u}{e}{w}=\tpc{u}{e}{w-v}+\tp{u}{e}{v}=0$. By Peirce calculus again, $\tp{e}{u}{w}=0$ as well.

We next show that $u\perp e$.   By Peirce calculus with
respect to $w$, $\tp{e}{u}{w}=0$. By the main identity,
$\tp{u}{e}{e}=\tp{u}{e}{\tp{e}{w}{w}}=\tp{\tp{u}{e}{e}}{w}{w}-\tp{e}{\tp{e}{u}{w}}{w}+\tp{e}{w}{\tp{u}{e}{w}}$.
The last two terms are zero and since $\tp{u}{e}{e}\in B_1(w)$, the first
term is equal to $\tp{u}{e}{e}/2$.   Hence $\tp{u}{e}{e}=0$ and $u\perp
e$.

Finally, we show that $e=0$.  Note first that by the Peirce calculus
$\tp{u}{v}{e}\in B_1(w)\cap B_1(u)\cap B_0(v)$ so $\tp{u}{v}{e}=0$ and by
Peirce calculus with respect to $w$, $\tp{v}{u}{e}=0$.
Hence, by the main identity,
$0=\tp{v}{u}{\tp{u}{v}{e}}=\tp{\tp{v}{u}{u}}{v}{e}-\tp{u}{\tp{u}{v}{v}}{e}+\tp{u}{v}{\tp{v}{u}{e}}=\tp{v}{v}{e}/2-\tp{u}{u}{e}/2+0=e/4$.

   From the fact just proved, namely, that $B_1(v)\cap B_1(w-v)=0$, it  
follows from
the joint Peirce decomposition that
$B_2(w)=B_2(v)\oplus B_2(w-v)$, which by \cite[Theorem 4.2]{HornMS}
implies that  $B=C\oplus D$ where $C$ and $D$ are orthogonal
weak$^*$-closed ideals and $u$ is a maximal tripotent of
$C=B_{2}(v)\oplus B_{1}(v)$.  It follows
from this and Lemma \ref{lem:5.1} that $u$ is a central tripotent of $B$. \qed

\medskip

   The proof  of the following remark is identical to the proofs of
Lemmas~\ref{lem:2} and ~\ref{lem:3}.  Recall that, as noted above, those
two lemmas are valid without assuming  the maximality of $w$ there and $u$ here.

\begin{remark}\label{prop:3.10}
Let $w$ be a maximal tripotent and let $u\in B_1(w)$ be a tripotent.
Assume that $u$ is not a central tripotent of $B$ and that $w_u\ne w$. Let
$a$ be a non-zero tripotent of $B_1(u)\cap B_0(w_u)\cap B_1(w)$ (which
is non-zero by Lemma~\ref{lem:horn}). Then
   $u_a\ (:=2\tp{a}{a}{u})\hbox{ is a
tripotent }\le u$ by Lemma \ref{lem:3},
noting that $\tp{a}{u}{a}=0$ by Peirce calculus with respect to $w_u$.  
Also $u_{a}$ lies in $B_{1}(a)$ by Peirce Calculus since $P_{2}(a)u=\tp{a}{\tp{a}{u}{a}}{a}=0$.
   \end{remark}

\begin{lemma}\label{lem:4.12}
With the notation of Remark~\ref{prop:3.10},
$w_{u_a}\top u_a$
\end{lemma}
\pf\
By assumption, $a\in B_1(w)$.  Therefore $u_a:=2\tp{u}{a}{a}\in B_1(w)$
and the result follows from Proposition~\ref{prop:3.2}.\qed

\begin{proposition}\label{prop:12261}
Let $B$ be a $JBW^*$-triple with no direct summand of the form
$L^\infty(\Omega,H)$ where $H$ is a Hilbert space of any positive
dimension.  Then every tripotent of $B$ is the supremum of the non-central
tripotents that it contains.
\end{proposition}
\pf\
Given a tripotent $u$ in $B$, let $v$ denote the supremum of all
non-central tripotents majorized by $u$, or zero, if there are none.  By
the definition of $v$,  $u-v$ is a central tripotent and any tripotent
majorized by $u-v$ is also a central tripotent. Hence $u-v$ is an abelian
tripotent, that is,
$B_2(u-v)$ is  associative and hence a commutative $C^*$-algebra.

   Thus
$B_2(u-v)\oplus B_1(u-v)$ is a weak$^*$-closed ideal containing a complete (=maximal) 
abelian tripotent, namely $u-v$.  By \cite[Theorem 2.8]{HornMZ},
$B_2(u-v)\oplus B_1(u-v)$ is a direct sum of spaces of the form
$L^\infty(\Omega_m,H_m)$ where $H_m$ is a Hilbert space of dimension $m$ for a family of
cardinal numbers $m$.  Since $u-v$ is a central tripotent, $B$ contains the weak*-closed ideal
$B_2(u-v)\oplus B_1(u-v)$ as an $\ell^\infty$-summand, contradicting our assumption.
Thus $u=v$.\qed

\subsection{Pullback of the Peirce 1-space}

We are now ready to prove the main result of this section.

\begin{theorem}\label{thm:2} 
Assume that $B$ has no direct summand of the form $L^\infty(\Omega,H)$
where $H$ is a Hilbert space of dimension at least two.  Suppose $w  
\in B$ is a maximal tripotent. Then
$u_\psi\in A_1(w_\psi)$ if $u\in B_1(w)$.
\end{theorem}
\pf\
Since commutative $JBW^*$-triples have no Peirce 1-spaces, it follows  
easily using a joint Peirce decomposition of $w$ that we may assume
$B$ also has no summands $L^\infty(\Omega)$, so that the hypothesis
of Proposition~\ref{prop:12261} holds. Thus we can  write $u=\sup_{\lambda\in\Lambda}
u_\lambda$ where each $u_\lambda$ is a non-central tripotent belonging to $B_1(w)$. Then by
Lemma \ref{lem:horn} and Remark~\ref{prop:3.10}, for each $\lambda\in\Lambda$, $v_\lambda:=\sup_a  
(u_\lambda)_a$ exists, where the supremum is over all non-zero tripotents $a$ in $B_1(u_\lambda)\cap B_0(w_{u_\lambda})\cap B_1(w)$.

We claim that $u=\sup_{\lambda\in\Lambda}
v_\lambda$                                                            .
Indeed, setting $v=\sup_\lambda
v_\lambda$, if $v\ne u$ we would have that  $u-v$ is the supremum of
non-central tripotents majorized by $u-v$ and hence by $u$.  Let
$u_{\lambda_0}$ be one of these non-central tripotents.  Then
$v_{\lambda_0}\le u_{\lambda_0}\le u-v$ which contradicts $v=\sup_\lambda
v_\lambda$. This proves the claim.

Explicitly,  we have proved
\[
u=\sup_\lambda\sup_{a_\lambda} \{(u_\lambda)_{a_\lambda}:a_\lambda\in B_1(w)\cap B_1(u_\lambda)\cap B_0(w_{u_\lambda})\}
\]
and this is the same as
\[
u=\sup\{(u_\lambda)_{a_\lambda}:\lambda\in\Lambda, a_\lambda\in  B_1(u_\lambda)\cap B_0(w_{u_\lambda})\cap B_1(w) \}.
\]

In the rest of this proof, we shall use the fact, just established, that
$u$ is the supremum of  a family of tripotents $v_a$ for certain $v \le
u$ and certain
tripotents $a \in B_1(v)\cap B_0(w_v)\cap B_1(w)$ where, by the argument at the end of Remark \ref{prop:3.10},  $v_{a}$ lies  
in $B_{1}(a)$.  Note that Lemma~\ref{6.7} will be used several times, as
indicated below.

We note first that  $w_{v_a},v_a\in B_2(w_v+a)$ and $v_a\in B_1(w_v)$.  Indeed,
  from $v_a\le v$ we have from Lemma~\ref{lem:7} that $w_{v_a}\le w_v$ so
$w_{v_a}\in B_2(w_v)\subset B_2(w_v+a)$.  On the other hand, by
Lemma~\ref{lem:4}, $v_a\in B_1(a)\cap B_2(v)\subset B_1(a)\cap
B_1(w_v)\subset B_2(w_v+a)$.

We claim next that $(v_a)_\psi\in A_1((w_u)_\psi)$.   Indeed,
since  by Lemma~\ref{lem:7},  $w_{v} \perp w_{u}-w_{v}$, we have by Remark~\ref{rem:4.19} and the joint Peirce decomposition,
\begin{equation}\label{eq:02021}
A_{1}((w_{v})_{\psi}) \cap A_{0}((w_{u}-w_{v})_{\psi})\subset A_1((w_u)_\psi).
\end{equation}
Since $w_v,v_a\in B_2(w_v+a)$ and
$\tp{w_v}{w_v}{v_a}=v_a/2$, it follows (using Lemma~\ref{6.7}) that
$\tpc{(w_v)_\psi}{(w_v)_\psi}{(v_a)_\psi}=(v_a)_\psi/2$ so $(v_a)_\psi$
lies in $A_{1}((w_{v})_{\psi})$.
Also, $v \perp w-w_{v}$ since
\begin{eqnarray*}
   \tpc{w-w_v}{w-w_v}{v}&=&
\tp{w}{w}{v}-\tp{w_v}{w}{v}-\tp{w}{w_v}{u}+\tp{w_v}{w_v}{v}\\
   &=& \tp{w}{w}{v}-\tp{w_v}{w_v}{v}-\tp{w_v}{w_v}{v}+\tp{w_v}{w_v}{v}\\
   &=&\tp{w}{w}{v}-\tp{w_v}{w_v}{v}=v/2-v/2=0.
   \end{eqnarray*}
Hence $v_{a}\le v$ lies in $A_{0}(w-w_{v}) \subseteq
A_{0}(w_{u}-w_{v})$. Embedding $v_{a}$ and $w_{u}-w_{v}$ in
$B_{2}(v_{a}+ w_{u}-w_{v})$, we see that $(v_a)_\psi$ lies in
$A_{0}((w_{u}-w_{v})_{\psi})$ and the  claim follows from (\ref{eq:02021}).

We now have from Lemma \ref{lem:4.2} and
Lemma~\ref{lem:3.13} that $u_\psi\in A_1((w_u)_\psi)$.
As before, $u\perp(w-w_u)$, so application of 
 Lemma~\ref{6.7}  and Remark~\ref{rem:4.19}  yields $u_\psi\in A_0((w-w_u)_\psi)$.
Finally, $u_\psi\in A_1((w_u)_\psi)\cap A_0((w-w_u)_\psi)\subset A_1(w_\psi)$.\qed

\section{The space of local multipliers}

We retain the notation of the previous two sections, that is,
$\psi:B_*\rightarrow A_*$ is a linear isometry, where $A$ and $B$ are
\jbwst s and $w$ is an extreme point of $B$  giving rise to the objects
$v,M,s$ in $A$. We also assume that $B$ satisifes the condition in
Theorem~\ref{thm:2}, that is, it
has no direct summand of the form $L^\infty(\Omega,H)$ where $H$ is a
Hilbert space of dimension at least two.

\begin{lemma}\label{lem:6.2}
$\psi[B_1(w)_*]\subset A_1(s)_*$.
\end{lemma}
\pf\
If $f\in B_1(w)_*$, then $v_ f\in B_1(w)$ and by Lemma~\ref{lem:6.1} and
Theorem~\ref{thm:2}, $v_ {\psi(f)}=(v_ f)_\psi\in A_1(s)$.

   To show that $\psi(f)\in A_1(s)_*$, let $g=\psi(f)$ and Peirce decompose
it with respect to $s$: $g=g_2+g_1+g_0$.  Since
$\pair{g_0}{A_0(s)}=\pair{g}{A_0(s)}=\pair{f}{\psi^*[A_0(s)]}=0$ we have
$g_0=0$.  It remains to show $g_2=0$. We may assume that $\|f\|=1$.

   Since $g=g_2+g_1$ and $v_g\in A_1(s)$, $g_1(v_g)=g(v_g)=1=\|g\|\ge \|g_1\|$ so that
$\|g_1\|=1$ and $g_1\in A_2(v_g)_*$.
   Since obviously $g\in A_2(v_g)_*$, we have $g_2\in A_2(v_g)_*$.   By
\cite[Lemma 1.1]{FriRus85bis}, we have  $\|\lambda
g_2+g_1\|=\|g_2+g_1\|=1$ for every complex $\lambda$ of modulus 1.  The
local argument given in
   \cite[Theorem 3.1]{AkeRus70} can be easily extended to apply to
$JBW^*$-algebras  to show that since $g_1$ is a complex extreme point of the
unit ball of the predual of the $JBW^*$-algebra $A_2(v_g)$, we must have
$g_2=0$.
    \qed

\begin{corollary} \label{NOP2}
     $\psi^*(A_2(s))\subset B_2(w)$
\end{corollary}
\pf\
If $x\in A_2(s)$ let $\psi^*(x)=y_2+y_1 $ be the Peirce decomposition of
$\psi^*(x)$ with respect to $w$. If $f\in B_1(w)_*$, then
$\pair{f}{y_1}=\pair{f}{\psi^*(x)-y_2}=\pair{f}{\psi^*(x)}=\pair{\psi(f)}{x}=0$
since $\psi(f)\in A_1(s)_*$ and $x\in A_2(s)$.  Thus $y_1=0$.\qed

\medskip

In view of this Corollary, we may improve the statement of
Lemma~\ref{lem:11141} by replacing $V$ by $\psi^*$   We restate this
improved lemma here.

\begin{lemma}\label{lem:11082}
Let $x\in A_2(s)$ be such that $0\le x\le s$ and $\psi^*(x)$ is a
projection in $B_2(w)$.   Then $x\in M_2(s)$.
\end{lemma}

The following is the announced improvement of Lemma~\ref{lem:01081}.
\begin{lemma}\label{lem:description}
Let $u$ be a tripotent in $B$.   Then
\begin{description}
\item[(a)]
$\psi^*(A_1(u_\psi))\subset B_1(u)+B_0(u)$
\item[(b)] $\psi^*(A_j(u_\psi))\subset B_j(u)$ for $j=0,2$
\end{description}
\end{lemma}
\pf\
Part (a) and the case $j=0$ of part (b) have been proved in
Lemma~\ref{lem:01081}.

To prove the case $j=2$ of (b), note first that by Lemma~\ref{lem:11082}
$u_\psi\in M_2(s)$.  (Recall that
   $u_\psi\le s\le v$ where $v$ is a maximal tripotent of $A$ with
$\psi^*(v)=w$ and $w$ is a maximal tripotent majorizing $u$.)

   If $x\in A_2(u_\psi)$, then
$x=\tp{u_\psi}{\tp{u_\psi}{x}{u_\psi}}{u_\psi}$ and by definition of
multiplier and using Corollary \ref{NOP2}, $\psi^*(u_\psi\circ  
c)=u\circ \psi^*(c)$ for all $c\in
A_2(v)$.  Using this and the general formula
$\tp{x}{y}{x}=2x\circ(x\circ y^*)-y^*\circ x^2$ we obtain
$\psi^*\tp{u_\psi}{x}{u_\psi}=\tpc{u}{\psi^*(x)}{u}$. For the same
reason,
$\psi^*(x)=\tpc{u}{\psi^*\tp{u_\psi}{x}{u_\psi}}{u}=\tp{u}{\tpc{u}{\psi^*(x)}{u}}{u}\in
B_2(u)$, proving the case $j=2$ of (b).
   \qed

\medskip

\begin{lemma}\label{6.6}
Suppose $x\in A$. If $\psi^*(x^{2n+1})=(\psi^*(x))^{2n+1}$ for all
positive integers $n$, then $x=(\psi^*(x))_{\psi}+q$, where $q \perp
(\psi^*(x))_{\psi}$.
\end{lemma}
\pf\
We may assume $||x||=1$.  Let $W(x)$ be the $JBW^*$-triple generated by
$x$.   By assumption and weak$^*$-continuity, $\psi^*$ restricts to an
isomorphism of $W(x)$ onto $W(\psi^*(x))$. For each closed subset $S$ of
$(0,1]$ if we let $u_{S}\in W(x)$ be the corresponding spectral tripotent
for $x$, then $\psi^*(u_{S})$ is the spectral tripotent $v_{S}$ of
$\psi^*(x)$ (or zero, if $S$ has no intersection with the spectrum of
$\psi^*(x)$).

Choose a maximal tripotent $w \ge r(\psi^*(x))$. If $\psi^*(u_{S})$ is not
zero, then by Lemma~\ref{6.85}, $u_{S}=(v_{S})_{\psi}+q_{S}$ where $q_{S}$
is a tripotent which is perpendicular to $(v_{S})_{\psi}$.

Now suppose $S \cap T=0$ and $u_{S}$ and $u_{T}$ are non-zero.     Then
$u_{T} \perp u_{S}$ and hence $(u_{T})_\psi$ is perpendicular to
$(v_{S})_{\psi}$ and $q_{S}$ (Remark~\ref{rem:4.19}). By symmetry, $u_{S}$
is perpendicular to $(v_{T})_{\psi}$ and $q_{T}$. A simple calculation of
$0=\{u_{S},u_{S},u_{T}\}$ shows that $q_{S} \perp q_{T}$.

It follows by approximation that $x=(\psi^*(x))_{\psi}+q$, where $q \perp
(r(\psi^*(x)))_{\psi}$. Indeed,
approximate $x$ as a norm limit of finite sums $y=\sum\lambda_i u_{S_i}$  
with the $S_i$ disjoint,
and $\sum u_{S_i} = r(x)=r(y)$.  Then $y=\sum\lambda_i u_{S_i}=\sum  
\lambda_i[(v_{S_i})_\psi+q_{S_i}]=(\sum\lambda_i
v_{S_i})_\psi+\sum q_{S_i}=(\psi^*(y))_\psi+q$  
where, since $q_{S_i}\perp
(v_{S_j})_\psi$ for all $i\ne j$, the element $q=\sum q_{S_i}$ is
orthogonal to  
$(r(\psi^*(y)))_\psi=(\psi^*(r(y)))_\psi=\sum(v_{S_i})_\psi$ and hence  
orthogonal to $(\psi^*(y))_{\psi}$. The result follows from continuity.
\qed

\medskip

Note that by the spectral theorem, Theorem~\ref{thm:2} is valid for arbitrary elements $x\in
B_1(w)$.   We now extend
Theorem~\ref{thm:2} to not necessarily maximal tripotents.

\begin{lemma}\label{lem:6.3} 
If $u$ is any tripotent of $B$ and if $x\in B_1(u)$, then $x_\psi\in
A_1(u_\psi)$.
\end{lemma}
\pf\
Consider first a tripotent $v \in B_1(u)$. Write
\[
v_\psi=P_2(u_\psi)v_\psi+P_1(u_\psi)v_\psi+P_0(u_\psi)v_\psi:=(v_\psi)_2+(v_\psi)_1+(v_\psi)_0
\]
and take $f\in B_1(u)_*$ with $f(v)=1=\|f\|$.  Then by
Lemma~\ref{lem:description}
\begin{eqnarray*}
1&=&f(v)=\psi(f)(v_\psi)=\psi(f)((v_\psi)_2+(v_\psi)_1+(v_\psi)_0)\\
&=&f(\psi^*[(v_\psi)_2]+\psi^*[(v_\psi)_1]+\psi^*[(v_\psi)_0])\\
&=&f[\psi^*[(v_\psi)_1]]=\psi(f)[(v_\psi)_1].
\end{eqnarray*}

Therefore by Lemma~\ref{lem:6.1} (recalling that $v_g$ denotes the support tripotent of the normal functional $g$), $(v_\psi)_1\ge v_{\psi(f)}=(v_f)_\psi$.
By  Lemma~\ref{lem:4.2}, and the fact, already used in Lemma~\ref{lem:4.1}, that every tripotent is the supremum of a family of support tripotents of normal functionals,
\begin{equation}\label{eq:0115081}(v_\psi)_1\ge \sup_f v_{\psi(f)}=\sup_f
(v_f)_\psi=(\sup_f v_f)_\psi=v_\psi=(v_\psi)_2+(v_\psi)_1+(v_\psi)_0.
\end{equation}

For notation's sake, let $y=v_\psi$.  The meaning of (\ref{eq:0115081}) is
that $(y_2+y_0)\perp y$, or
$D(y_2+y_0,y_2+y_0)(y_2+y_1+y_0)=0$.
   This yields, upon expansion and comparison of Peirce components, that
$\tp{y_2}{y_2}{y_2}=0=\tp{y_0}{y_0}{y_0}$ so that $y_2=y_0=0$.
Thus, $v_{\psi}$ lies in $A_1(u_\psi)$.

The lemma follows easily for an arbitrary $x\in B_1(u)$ by considering the
spectral decomposition of $x$.
\qed

\begin{lemma}\label{lem:6.4}
Let $u$ and $v$ be compatible tripotents in $B$ (in particular, if $u$ is a tripotent in $B_{1}(v)$) and let
$x$ be an element in $B_{2}(v)$.  Then
\[
P_j(u_\psi)x_\psi=(P_j(u)x)_\psi\hbox{ for }j=0,1,2.
\]
In particular  $P_j(u_\psi)x_\psi$ is a multiplier for $j=0,1,2$.
\end{lemma}
\pf\
Since $u$ and $v$ are compatible, $P_j(u)x=P_2(v)P_j(u)x\in B_2(v)$ so that by Lemma~\ref{6.7}, 
 \begin{equation}\label{eq:02022}
x_\psi=(P_2(u)x+P_1(u)x+P_0(u)x)_\psi=(P_2(u)x)_\psi+(P_1(u)x)_\psi+(P_0(u)x)_\psi.
\end{equation}

From Lemma~\ref{lem:6.3}, $(P_1(u)x)_\psi\in A_1(u_\psi)$ and by Remark~\ref{rem:4.19},
$(P_0(u)x)_\psi\in A_0(u_\psi)$.  Again by Lemma~\ref{6.7}, $$(P_2(u)x)_\psi=(\tp{u}{\tp{u}{x}{u}}{u})_\psi=(\tp{u}{\tpc{u}{P_2(u)x}{u}}{u})_\psi=(\tp{u_\psi}{\tpc{u_\psi}{(P_2(u)x)_\psi}{u_\psi}}{u_\psi}),$$
so that $(P_2(u)x)_\psi\in A_2(u_\psi)$.

 By the uniqueness of Peirce decompositions and (\ref{eq:02022}),
$P_j(u_\psi)x_\psi=(P_j(u)x)_\psi$.\qed

\section{Proof of the main result}

We again assume in this section that $B$ satisifes the condition in
Theorem~\ref{thm:2}.

\begin{lemma}\label{6.5}
Suppose $v$ is a tripotent in $B$.
Further suppose that $x$ is a tripotent  in $B_{1}(v)$ with $\{x,v,x\}=0$
and $\{x_\psi,v_\psi,x_\psi\}=0$.  Then
$\psi^*\{x_{\psi},x_{\psi},v_{\psi}\}=\{x,x,v\}$. Furthermore,
$\{x_{\psi},x_{\psi},v_{\psi}\}=y_{\psi}$ for some $y \in B$.
\end{lemma}

\pf\
We note first that, as shown in Lemma~\ref{lem:3}, $p:=2D(x,x)v$ is a
self-adjoint projection in $B_{2}(v)$.
By Peirce arithmetic, using the assumption $\tp{x}{v}{x}=0$,
$p$ lies in $B_{1}(x)$ and by Lemma \ref{lem:6.3}, $p_{\psi}$
lies in $A_{1}(x_{\psi})$. By this fact, the compatibility of $p_\psi$ and $x_\psi$, and the fact that $p_\psi\le v_\psi$, we have
\[
2D(p_{\psi},p_\psi)D(x_{\psi},x_\psi)v_{\psi}=2D(x_{\psi},x_\psi)D(p_{\psi},p_\psi)v_{\psi}=2D(x_{\psi},x_\psi)p_{\psi}=p_{\psi}.
\]

Similarly to the calculation above,  $q:=2\{x_{\psi},x_{\psi},v_{\psi}\}$
is a self-adjoint projection in $A_{2}(v_{\psi})$ and since $q\circ
p_\psi=2\tp{\tp{x_\psi}{x_\psi}{v_\psi}}{v_\psi}{p_\psi}=2D(p_\psi,p_\psi)D(x_\psi,x_\psi)v_\psi=p_\psi$,
$q\ge p_{\psi}$ and it follows that
$\psi^* (q) \ge p$. 

Now
$D(x,x)(v-p)=\{xxv\}-\{xxp\}=p/2-p/2=0$. Hence, $x_{\psi}$ is orthogonal
to $v_{\psi}-p_{\psi}$.
By this orthogonality and compatibility, and since $p_\psi\le v_\psi\le
w_\psi$ ($w$ is a maximal tripotent majorizing $v$) so that
$\tp{p_\psi}{p_\psi}{v_\psi}=p_\psi$,
\begin{eqnarray*}
D(v_{\psi}-p_{\psi},v_{\psi}-p_{\psi})D(x_{\psi},x_{\psi})v_{\psi}&=&D(x_{\psi},x_{\psi})D(v_{\psi}-p_{\psi},v_{\psi}-p_{\psi})v_{\psi}\\
&=&D(x_{\psi},x_{\psi})(v_{\psi}-p_{\psi})=0,
\end{eqnarray*}
showing $v_{\psi}-p_{\psi}$ is orthogonal to
$q$.  We then have
$\|v-p\pm\psi^*(q)\|\le \|v_\psi-p_\psi\pm q\|=1$ so that  $v-p$ is
orthogonal to $\psi^*(q)$.  Since, as shown above,  $\psi^*(q)\ge  p$,  it follows (using
Lemma~\ref{lem:description} to ensure that $\psi^*(q)\in B_2(v)$) that
$\psi^*(q)= p$. This proves the first
statement. The second follows immediately from Lemma~\ref{lem:11082} since
$v_{\psi}$ is majorized by $w_{\psi}$ for a maximal tripotent $w \in B$ and
$\psi^*$ takes the positive element $2\{x_{\psi},x_{\psi},v_{\psi}\} \in
A_{2}(w_{\psi})$ to a
projection in $B_{2}(w)$.  \qed

\begin{lemma}\label{6.8}
Suppose that $y$ and $z$ lie in $B_{2}(w)$ for a maximal tripotent $w$ and
that $x$ lies in $B_{1}(w)$. Then $\{x_{\psi},y_{\psi},z_{\psi}\}$ is a
multiplier in $A_{1}(w_{\psi}) \cap A_{2}([r(x)+r(z_{0})]_{\psi})$ (where $z_0=P_0(r(x))z$), and
$\psi^{\ast}\{x_{\psi},y_{\psi},z_{\psi}\}=\{x,y,z\}$.
\end{lemma}
\pf\
Suppose first that $x$ is a tripotent.  Let  $y_j$ denote $P_j(x)y$ and
$(y_\psi)_j=P_j(x_\psi)y_\psi$ for
$j=0,1,2$. Similarly for $z$.
By Lemma~\ref{lem:6.4}, replacing $u,v,x$ there by $x,w,y$ respectively,
we have in particular that  $(y_{1})_{\psi}=(y_{\psi})_{1}$ and similarly
$(z_{1})_{\psi}=(z_{\psi})_{1}$.

Note that in the expansion

\[
\{x_\psi, y_\psi, z_\psi\}=\{x_\psi,\sum_i
(y_\psi)_i,\sum_j(z_\psi)_j\}=\sum_{i,j}\{x_\psi,(y_\psi)_i,(z_\psi)_j\},
\]
seven of the nine terms are zero, five of them since
$y_2=\tp{x}{\tpc{x}{y}{x}}{x}=0$ by the maximality of $w$ (so also
$z_2=0$), and two others since $x_\psi\perp (y_\psi)_0$.
   Hence
\begin{equation}\label{eq:12221}
\{x_{\psi},y_{\psi},z_{\psi}\}=\{x_{\psi},(y_{1})_{\psi},(z_{1})_{\psi}\}+\{x_{\psi},(y_{1})_{\psi},(z_{0})_{\psi}\}.
\end{equation}

Let $u_{S}$ be a spectral tripotent of $y_{1}$. By Peirce calculus with
respect to $w$ and $w_{\psi}$,
$\{u_{S},x,u_{S}\}=0$ and $\{(u_{S})_{\psi},x_{\psi},(u_{S})_{\psi}\}=0$.
Therefore, by Lemma~\ref{6.5},
   $\{x_{\psi},(u_{S})_{\psi},(u_{S})_{\psi}\}$ is a multiplier in
$A_2(x_\psi)$ and
$\psi^*\{x_{\psi},(u_{S})_{\psi},(u_{S})_{\psi}\}=\tp{x}{u_S}{u_S}$.
Passing to the limit using the spectral theorem shows that
$\{x_{\psi},(y_{1})_{\psi},(y_{1})_{\psi}\}$ is a multiplier in
$A_2(x_\psi)$ and
$\psi^*\{x_{\psi},(y_{1})_{\psi},(y_{1})_{\psi}\}=\tp{x}{y_1}{y_1}$. Of
course, the same holds for $z$:
$\{x_{\psi},(z_{1})_{\psi},(z_{1})_{\psi}\}$ is a multiplier in
$A_2(x_\psi)$ and
$\psi^*\{x_{\psi},(z_{1})_{\psi},(z_{1})_{\psi}\}=\tp{x}{z_1}{z_1}$.

By  Lemma~\ref{6.7}, $(y_{1})_{\psi}+(z_{1})_{\psi}=(y_{1}+z_{1})_{\psi}$.
   Hence the same statement
holds for
$\{x_{\psi},(y_{1})_{\psi}+(z_{1})_{\psi},(y_{1})_{\psi}+(z_{1})_{\psi}\}$.
   Thus the statement
holds for
$\{x_{\psi},(y_{1})_{\psi},(z_{1})_{\psi}\}+\{x_{\psi},(z_{1})_{\psi},(y_{1})_{\psi}\}$.
   Explicitly,
$\{x_{\psi},(y_{1})_{\psi},(z_{1})_{\psi}\}+\{x_{\psi},(z_{1})_{\psi},(y_{1})_{\psi}\}$
is a multiplier
   in $A_2(x_\psi)$ and
    $$\psi^*(\{x_{\psi},(y_{1})_{\psi},(z_{1})_{\psi}\}+\{x_{\psi},(z_{1})_{\psi},(y_{1})_{\psi}\}
)=\tp{x}{y_1}{z_1}+\tp{x}{z_1}{y_1}.
$$

Replacing $z$ by $iz$ shows  that the statement holds for
$\{x_{\psi},(y_{1})_{\psi},(z_{1})_{\psi}\}$ and
$\{x_{\psi},(z_{1})_{\psi},(y_{1})_{\psi}\}$ individually.  This proves,
in the case that $x$ is a tripotent,  that the first  term in
(\ref{eq:12221}) is a multiplier in $A_{2}(x_{\psi}) \cap
A_{1}(w_{\psi})$ and and $\psi^*$ is multiplicative on
this term.

We now consider the second term in (\ref{eq:12221}), still in the case
that $x$ is a tripotent.
Since $x\perp z_0$ (recall that $z_0=P_0(x)z$), we can choose a maximal
tripotent $w'$ such that
$B_2(x+r(z_0))\subset B_2(w')$, so that $x_\psi$ and $(z_0)_\psi$ are
multipliers in $A_2(x_\psi+r(z_0)_\psi)=A_2([x+r(z_0)]_\psi)\subset
A_2(w'_\psi)$. We next note that for every $a\in A$,
\begin{equation}\label{eq:12222}
\psi^*\{x_\psi,a,(z_0)_\psi\}=\{x,\psi^*(a),z_0\}.
\end{equation}
Indeed, by Peirce calculus
$\{x_\psi,a,(z_0)_\psi\}=\psi^*\{x_\psi,P_2(w'_\psi)a,(z_0)_\psi\}$ and by
properties of multipliers and the Jordan algebra relation
\begin{equation}\label{eq:12223}
\tp{a}{b}{c}=(a\circ b^*)\circ c+(c\circ b^*)\circ a-(a\circ c)\circ b^*,
\end{equation}  (cf.\ Lemma~\ref{lem:2.2bis}), and
Lemma~\ref{lem:description},
\begin{eqnarray*}
\psi^*\{x_\psi,a,(z_0)_\psi\}&=&\psi^*\{x_\psi,P_2(w'_\psi)a,(z_0)_\psi\}\\
&=&\{x,\psi^*(P_2(w'_\psi)a),(z_0)_\psi\}\\
&=&
\{
x,P_2(w')\psi^*(a),z_0\}\\
&=&\{x,\psi^*(a),z_0\}.
\end{eqnarray*}
In particular, $\psi^*\{x_\psi,(y_1)_\psi,(z_0)_\psi\}=\{x,y_1,z_0\}$ so
that
$$\psi^*\{x_\psi,y_\psi,z_\psi\}=\{x,y_1,z_1\}+\{x,y_1,z_0\}=\{x,y,z\}.$$

It remains to show that $\{x_\psi,(y_1)_\psi,(z_0)_\psi\}$ is a multiplier.
By the joint Peirce decomposition and the relation
$D(u,u)=P_2(u)+P_1(u)/2$,
\begin{eqnarray*}
P_{2}(x_{\psi}+r(z_{0})_{\psi})(y_{1})_{\psi}
&=&[P_{2}(x_{\psi})+P_2(r(z_{0})_{\psi})+ P_1(x_\psi)P_1(r(z_0)_\psi)
](y_{1})_{\psi}\\
&=&P_1(r(z_0)_\psi)  (y_{1})_{\psi}\\
&=&[2D(r(z_0)_\psi,r(z_0)_\psi)-2P_2(r(z_0)_\psi)](y_1)_\psi\\
&=&2D(r(z_0)_\psi,r(z_0)_\psi)(y_1)_\psi.
\end{eqnarray*}

   The right side of the preceding equation is a  triple product  of
multipliers in $A_{2}(w_{\psi})$ and is hence a
   multiplier in $A_{2}(w_{\psi})$ by (\ref{eq:12223}) and the fact that the
multipliers form a Jordan algebra.
   Hence $P_{2}(x_{\psi}+r(z_{0})_{\psi})(y_{1})_{\psi}$ is a multiplier in
$A_2(w_\psi)$.
      Since
$\tp{x_\psi}{(y_1)_\psi}{(z_0)_\psi}=2P_{2}(x_{\psi}+r(z_{0})_{\psi})(y_{1})_{\psi}$,
   using Lemma \ref{6.7}, $\tp{x_\psi}{(y_1)_\psi}{(z_0)_\psi}$ is a
multiplier in $A_{2}([x+r(z_{0})]_{\psi})$.

Now let $x$ be an
arbitrary element of $B_1(w)$.  Approximate it by  sums
$\tilde x=\sum
\lambda_iu_i$ where the elements $u_i\in B_1(w)$ are orthogonal spectral
tripotents with $\sum u_i = r(x)$.  Decomposing $y$ and $z$ with respect
to $r(x)=r(\tilde x)$, it follows as in (\ref{eq:12221}) that
\begin{equation}\label{999}
\tp{\tilde x_\psi}{y_\psi}{z_\psi}=\{\tilde
x_{\psi},(y_{1})_{\psi},(z_{1})_{\psi}\}+\{\tilde
x_{\psi},(y_{1})_{\psi},(z_{0})_{\psi}\}.
\end{equation}

By the previous discussion,
$\{(u_{i})_{\psi},(y_{1})_{\psi},(z_{1})_{\psi}\}$, which lies in
$A_{2}(r(x)_{\psi})$ by Peirce calculus, is a sum of a multiplier in
$A_{2}((u_i)_{\psi}) \subseteq A_{2}(r(x)_{\psi})$ and a multiplier in
$A_{1}(w_{\psi})$ which must thus also lie in $A_{2}(r(x)_{\psi})$.
Also, $\psi^{\ast}$ is multiplicative on these products.  Hence the
first term in (\ref{999}) is a multiplier in $A_{2}(r(x)_{\psi})
\subseteq A_{2}([r(x)+r(z_{0})]_{\psi})$ and $\psi^{\ast}$ is
multiplicative on it.

The second term equals $\sum\lambda_i
\{(u_{i})_{\psi},(y_{1})_{\psi},(z_{0})_{\psi}\}$. Since $z_{0} \perp
u_{i}$ the same argument used above shows that
$\{(u_{i})_{\psi},(y_{1})_{\psi},(z_{0})_{\psi}\}$ is a multiplier in
$A_{2}([u_{i}+r(z_{0})]_{\psi}) \subseteq A_{2}([r(x)+r(z_{0})]_{\psi})$
and  that $\psi^{\ast}$ is multiplicative on these products. Hence
the second term in (\ref{999}) is a multiplier in $A_{2}(r(x)_{\psi})
\subseteq A_{2}([r(x)+r(z_{0})]_{\psi})$ and $\psi^{\ast}$ is
multiplicative on it. The lemma follows.
\qed

\begin{lemma}\label{6.88}
If $q$ lies in $A_{0}(v_{\psi})$ for some maximal tripotent $v \in B$,
then $\psi^*\{q,q,x\}=0$ and $\psi^*\{q,x,y\}=0$ for all $x,y\in A$. Also,
$q \perp x_{\psi}$ for all $x \in B$, that is, $A_{0}(v_{\psi})\perp
\{x_\psi:x\in B\}$.
\end{lemma}
\pf\
Let $z$ be a maximal tripotent in $A_0(v_\psi)$ such that $q/\|q\|$ is a
self-adjoint element with respect to $z$ (see \cite[Lemma
3.12(1)]{HornMS}).  Clearly $v_{\psi}+z$ is maximal. Because $\psi^*$  
preserves orthogonality with $v_\psi$
and $v$ is maximal, $\psi^*(q)=\psi^*(z)=0$ and therefore $\psi^*$ maps the
self-adjoint element $v_\psi+q/\|q\|$ to the unit $v$ of $B_2(v)$ and  
maps $v_{\psi}+z$ to $v$. By
Corollary~\ref{cor:2.2}, $v_\psi+q/\|q\|$ is a multiplier in
$A_2(v_\psi+z)$.  Since $v_\psi$ is  a multiplier there, so is $q$.  On
the other hand, if we let $x=x_2+x_1+x_0$ be its Peirce decomposition with
respect to $v_\psi$, then $\tp{q}{q}{x}=\tpc{q}{q}{x_1+x_0}$ so that
$\psi^*\tp{q}{q}{x}=\psi^*\tp{q}{q}{x_1}$ since $\tp{q}{q}{x_0}\in
A_0(v_\psi)$.  If we now expand $x_1$ in its Peirce decomposition with
respect to $z$, say $x_1=(x_1)_2+(x_1)_1+(x_1)_0$, then
$\tp{q}{q}{x_1}=\tpc{q}{q}{(x_1)_2+(x_1)_1}$ and since $v_\psi$ and $z$
are compatible,
$(x_1)_2+(x_1)_1\in A_2(z)+A_1(z)\cap A_1(v_\psi)\subset A_2(v_\psi+z)$.
Since $q$
is a multiplier in $A_2(v_\psi+z)$, we now have
$\psi^*\tp{q}{q}{x_1}=\tpc{\psi^*(q)}{\psi^*(q)}{\psi^*((x_1)_2+(x_1)_1)}=0$,
proving that $\psi^*\tp{q}{q}{x}=0$.

Letting $x,y\in A$ and Peirce decomposing them with respect to $v_\psi$,
we have
\begin{equation}\label{eq:12224}
\psi^*\tp{q}{x}{y}=\psi^*\tpc{q}{x_1+x_0}{y_2+y_1+y_0}=\psi^*\tpc{q}{x_0}{y_1+y_0}+\psi^*\tp{q}{x_1}{y_2}.
\end{equation}

Since $\tp{q}{x_1}{y_2}\in A_1(z)$ (by Peirce calculus), we have
$\tp{q}{x_1}{y_2}=2\tpc{z}{z}{\tp{q}{x_1}{y_2}}=2\tpc{z}{v_\psi+
z}{\tp{q}{x_1}{y_2}}$ and therefore, since $z$ is a multiplier in
$A_2(v_\psi+z)$, $\psi^*\tp{q}{x_1}{y_2}=\psi^*(z)\circ
\psi^*\tp{q}{x_1}{y_2}=0$. Thus the second term on the right side of
(\ref{eq:12224}) is zero.

For the first term on the right side of (\ref{eq:12224}), we have
\begin{equation}\label{eq:12231}
\psi^*\tpc{q}{x_0}{y_1+y_0}=\psi^*\tpc{q}{x_0}{y_1}+\psi^*\tpc{q}{x_0}{y_0}
\end{equation}
and the second term in (\ref{eq:12231}) is zero since
$\tpc{q}{x_0}{y_0}\in A_0(v_\psi)$.   Peirce decomposing $x_0$ and $y_1$
with respect to $z$ and expanding the first term in (\ref{eq:12231}) leads
to
\begin{eqnarray*}
\psi^*\tpc{q}{x_0}{y_1}&=&\psi^*\tpc{q}{(x_0)_2}{(y_1)_2}+\psi^*\tpc{q}{(x_0)_2}{(y_1)_1}\\
&+&\psi^*\tpc{q}{(x_0)_1}{(y_1)_1}+\psi^*\tpc{q}{(x_0)_1}{(y_1)_0}.
\end{eqnarray*}

The first and third terms here are zero since $(y_1)_2$ and
$\tpc{q}{(x_0)_1}{(y_1)_1}$ belong to $A_1(v_\psi)\cap A_2(z)$, which is
zero since $v_\psi\perp z$.
The second term is zero since $\tpc{q}{(x_0)_2}{(y_1)_1}$ lies in  
$A_{1}(v_{\psi}) \cap A_{1}(z) \subseteq A_{2}(v_{\psi}+z)$ and
$\tpc{q}{(x_0)_2}{(y_1)_1}=2\tpc{z}{z}{\tpc{q}{(x_0)_2}{(y_1)_1}}=2\tpc{z}{v_\psi+z}{\tpc{q}{(x_0)_2}{(y_1)_1}}$
so that $\psi^*\tpc{q}{(x_0)_2}{(y_1)_1}=\psi^*(z)\circ
\psi^*\tpc{q}{(x_0)_2}{(y_1)_1}=0$.  The proof that the fourth term is  
zero is similar. This proves that
$\psi^*\tp{q}{x}{y}=0$.

To prove the last statement, it may be assumed that both $q$ and $x$ are
tripotents. Decompose $x_\psi$ with respect to $q$:
$x_\psi=(x_\psi)_2+(x_\psi)_1+(x_\psi)_0$ and note that by the first two
parts of this lemma, $\psi^*((x_\psi)_2+ (x_\psi)_1)=0$, so that
$\psi^*((x_\psi)_0) =x$.
By Lemma~\ref{6.85}, $(x_\psi)_0=x_\psi+\tilde q$ where $\tilde q\perp
x_\psi$.  Thus $\tilde q=-(x_\psi)_2-(x_\psi)_1$ is orthogonal to
$(x_\psi)_2+(x_\psi)_1+(x_\psi)_0$. Considering the components of  
$0=D((x_{\psi})_{2}+(x_{\psi})_{1},(x_{\psi})_{2}+(x_{\psi})_{1}+(x_{\psi})_{0})(x_{\psi})_{2}$
we immediately see that $(x_{\psi})_{2} \perp (x_{\psi})_{1}$ and  
$((x_{\psi})_{2})^{3}=0=(x_{\psi})_{2}$. Considering  
$0=D((x_{\psi})_{1},(x_{\psi})_{1}+(x_{\psi})_{0})(x_{\psi})_{1}$ we  
see that $(x_{\psi})_{1}=0$.
The lemma follows.
\qed

\begin{corollary}\label{6.87}
If $x \in B_{2}(w)$ for a maximal tripotent $w$ and $y,z \in B_{1}(w)$,
then $\{y_{\psi},x_{\psi},z_{\psi}\}=0$.
\end{corollary}
\pf\    Let $\alpha:=\{y_{\psi},x_{\psi},z_{\psi}\}$.  By Peirce calculus
with respect to $w_\psi$, $\alpha\in A_0(w_\psi)$ so by
Lemma~\ref{6.88}, $\{y_\psi,z_\psi,x_\psi\}\perp \alpha$.
   By the main identity,
$$\tp{\alpha}{\alpha}{\alpha}=\mi{\alpha}{\alpha}{y_\psi}{x_\psi}{z_\psi}
$$
and each term is zero, hence $\alpha=0$.\qed

\begin{lemma}\label{6.9}
Suppose $x_{\psi}$ is a multiplier in
$A_{1}(w_{\psi})$ for a maximal tripotent $w \in B$ and that $y_{\psi}$ is
a multiplier in $A_{2}(w_{\psi})$. Then $\{x_{\psi},x_{\psi},y_{\psi}\}$
is a multiplier and $\psi^*$ is multiplicative on this product.
\end{lemma}
\pf\
Suppose first that $x$ is a tripotent.   By Corollary~\ref{6.87},
$\tp{x_\psi}{y_\psi}{x_\psi}=0$ and hence
$P_2(x_\psi)y_\psi=0$.   Then by Lemma~\ref{lem:6.4},
\begin{eqnarray*}
   \tp{x_\psi}{x_\psi}{y_\psi}&=&D(x_\psi,x_\psi)y_\psi\\
   &=&(P_2(x_\psi)+P_1(x_\psi)/2)y_\psi\\
   &=&P_1(x_\psi)y_\psi/2=(P_1(x)y)_\psi/2,
\end{eqnarray*}
proving that $\{x_{\psi},x_{\psi},y_{\psi}\}$ is a multiplier.
Moreover,
$\psi^*\tp{x_\psi}{x_\psi}{y_\psi}=P_1(x)y/2=(2D(x,x)-2P_2(x))y/2=\tp{x}{x}{y}$,
since by Peirce calculus with respect to the maximal tripotent w,
$\{xyx\}=0$.

For the general case it suffices to assume that $x$ is a finite sum
$\sum\lambda_i x_i$ of pairwise orthogonal tripotents $x_i$ in $B_1(w)$.
By the special case just proved, $\tp{(x_i)_\psi}{(x_i)_\psi}{y_\psi}$ is
a multiplier and $\psi^*$ is multiplicative on it.  Therefore,
$$\{x_{\psi},x_{\psi},y_{\psi}\}=\sum\lambda_i^2
\tp{(x_i)_\psi}{(x_i)_\psi}{y_\psi}$$ is also a multiplier
and $\psi^*$ is multiplicative on it. \qed

\begin{lemma}\label{6.10}
Suppose that $z$ is a tripotent in $B$ and that $w$ is maximal tripotent
in $B$. Then, letting $z_{2}=P_{2}(w)z$ and $z_{1}=P_{1}(w)z$, we have
$z_{\psi}=(z_{2})_{\psi}+(z_{1})_{\psi}$
\end{lemma}
\pf\
It follows from Corollary~\ref{6.87} and Lemmas \ref{6.8} and \ref{6.9}
that $\psi^*[((z_{2})_{\psi}+(z_{1})_{\psi})^{3}]=z$.
Indeed,
$$((z_{2})_{\psi}+(z_{1})_{\psi})^{3}=\sum_{i,j,k=0}^1\tpc{(z_i)_\psi}{(z_j)_\psi}{(z_k)_\psi},$$
and $\psi^*$ is multiplicative on each term on the right side as follows.
For the terms corresponding to $(i,j,k)=(2,2,2)$ and $(1,1,1)$, this is
because $\psi^*$ is a Jordan homomorphism on the set of local multipliers.
For the terms corresponding to $(i,j,k)=(2,2,1)$ and $(1,2,2)$ (which are
the same), this is because of Lemma~\ref{6.8}.
For the terms corresponding to $(i,j,k)=(2,1,1)$ and $(1,1,2)$ (which are
the same), this is because of Lemma~\ref{6.9}.
For the term corresponding to $(1,2,1)$, this is because of
Corollary~\ref{6.87} and the maximality of $w$.
For the term corresponding to $(2,1,2)$, this is because of Peirce calculus.
Thus
$$\psi^*[((z_{2})_{\psi}+(z_{1})_{\psi})^{3}]=\sum_{i,j,k=0}^1\tp{z_i}{z_j}{z_k}=(z_2+z_1)^3=z^3=z,$$
as required.

Now if we Peirce decompose
$((z_{2})_{\psi}+(z_{1})_{\psi})^{3}$
with respect to $w_\psi$ we obtain
\begin{equation}\label{eq:02031}
P_2(w_\psi)[((z_{2})_{\psi}+(z_{1})_{\psi})^{3}]=((z_2)_\psi)^3+2\tpc{(z_2)_\psi}{(z_1)_\psi}{(z_1)_\psi},
\end{equation}
\begin{equation}\label{eq:02032}
P_1(w_\psi)[((z_{2})_{\psi}+(z_{1})_{\psi})^{3}]=((z_1)_\psi)^3+2\tpc{(z_2)_\psi}{(z_2)_\psi}{(z_1)_\psi},
\end{equation}
and
\[P_0(w_\psi)[((z_{2})_{\psi}+(z_{1})_{\psi})^{3}]=0.
\]

By Lemma~\ref{6.9}, the right side of (\ref{eq:02031})  is a sum of three
multipliers, and hence a multiplier itself in $A_{2}(w_{\psi})$.

On the other hand,
the first term on the right side of (\ref{eq:02032})  is obviously a multiplier in
$A_{2}(r(z_{1})_{\psi})\subseteq  
A_{2}([r(z_{1})+r(P_{0}(r(z_{1}))z_{2})]_{\psi})$.
By Lemma \ref{6.8}, the second
term is also a  multiplier in  
$A_{2}([r(z_{1})+r(P_{0}(r(z_{1}))z_{2})]_{\psi})$.
Hence the sum
is a multiplier.  It follows that $((z_{2})_{\psi}+(z_{1})_{\psi})^{3}$
is again a
sum of two multipliers $(z^{\prime}_{2})_{\psi}+(z^{\prime}_{1})_{\psi}$.
Repeating this argument, we see that $\psi^*[   (   (z_{2})_{\psi}
+(z_{1})_{\psi}     )^{3^{n}}]=z^{3^{n}}=z$, for  every $n$. Since
$C(x)=C(x^{3})$, we may use Lemma~\ref{6.6} to see that
$(z_{2})_{\psi}+(z_{1})_{\psi}=z_{\psi}+q$, where $q \perp z_{\psi}$ and
$\psi^*(q)=0$.

To show that $q=0$, suppose first that $z$ is maximal. It follows from
Lemma \ref{6.88} that $q \perp [(z_{2})_{\psi}+(z_{1})_{\psi}]$, from which
it follows that $q^{3}=0$, and $q=0$. Now suppose $z$ is a general
tripotent less than a maximal tripotent $v$. Let $u=v-z$. Then
$(z_{2})_{\psi}+(z_{1})_{\psi}+(u_{2})_{\psi}+(u_{1})_{\psi}=z_{\psi}+q
+u_{\psi}+p=v_{\psi}+p+q=(v_2)_\psi+(v_1)_\psi+p+q$.

Note that $(z_2)_\psi+(u_2)_\psi=(z_2+u_2)_\psi=(v_2)_\psi$ and therefore
\[
(v_2)_\psi+(z_1)_\psi+(u_1)_\psi=(v_2)_\psi+(v_1)_\psi+p+q
\]
which tells us that $p+q\in A_1(w_\psi)$.   Repeating this argument with
$-u$ instead of $u$ shows that
$p-q\in A_1(w_\psi)$ so that both $p$ and $q$ belong to $A_1(w_\psi)$.

   From $(z_2)_\psi+(z_1)_\psi=z_\psi+q$ with $q\in A_0(z_\psi)\cap
A_1(w_\psi)$ and $z_\psi=(z_\psi)_2+(z_\psi)_1$ we have $q\perp
(z_\psi)_1$; indeed,
$0=\tp{z_\psi}{q}{q}=\tp{(z_\psi)_2}{q}{q}+\tp{(z_\psi)_1}{q}{q}$ and both
terms are zero by Peirce calculus.

Thus $(z_1)_\psi=(z_\psi)_1+q$ with $q\perp (z_\psi)_1$ and therefore
\begin{equation}\label{eq:530}
r(z_1)_\psi=r((z_\psi)_1)+r(q)\hbox{ with }r(q)\perp r((z_\psi)_1).
\end{equation}

   From $q\perp z_\psi$ we have $r(q)\perp z_\psi$ and therefore
$\psi^*(r(q))\perp z$.  Since $\psi^*(r(q))$ lies in $A_{1}(w)$ by Lemma
\ref{lem:6.4},
a simple calculation as above shows $\psi^*(r(q))\perp z_2$ and
$\psi^*(r(q))\perp z_1$.  Finally, from
(\ref{eq:530}), $r(q)\le r(z_1)_\psi$ so that $r(q)\in A_2(r(z_1)_\psi)$
and $\psi^*(r(q))\in B_2(r(z_1))$.   But we already know that
$\psi^*(r(q))\in B_0(r(z_1))$, proving that
$\psi^*(r(q))=0$.

Now again by (\ref{eq:530}), $\psi^*(r((z_\psi)_1))=r(z_1)$ showing by
Lemmas \ref{lem:11141} and \ref{lem:2.2bis} that
$r(z_1)_\psi=r((z_\psi)_1)$, that is $r(q)=0$ and $q=0$.
\qed

\begin{theorem}\label{thm:3}
Let $\psi$ denote an isometry of $B_{\ast}$ into $A_{\ast}$
where $A$ and $B$ are JBW*-triples. Assume that $B$ has no
$L^{\infty}(\Omega,H)$ summand, where $H$ is a Hilbert space of dimension at least two. Let $C$ be the weak*-closure of the linear span of all multipliers: $C:=\overline{\hbox{sp}}^{w*}\{x_\psi|x
\in B\}$.  Then $C$  is a JBW*-subtriple of $A$, and $\psi^*$ restricted to
$C$ is a weak* bi-continuous isomorphism onto $B$ with inverse $x\mapsto x_\psi$ for
 $x\in B$.
\end{theorem}

\pf\
We first consider three tripotents $u,v$ and $w$ in $B$ and show that
$\{u_{\psi},v_{\psi},w_{\psi}\}$ is a sum of multipliers and that
$\psi^*$ is multiplicative on this product. Choose a maximal tripotent
$z \ge v$ and  decompose with respect to it: $u=u_{2}+u_{1}$ and
$w=w_{2}+w_{1}$.  It
follows from Lemma~\ref{6.10} and Corollary~\ref{6.87}, that the above
product equals
\[
\{(u_{2})_{\psi},v_{\psi},(w_{2})_{\psi}\}+\{(u_{1})_{\psi},v_{\psi},(w_{2})_{\psi}\}+\{(u_{2})_{\psi},v_{\psi},(w_{1})_{\psi}\}.
\]

The first product satisfies the desired conditions by the work in section 3.
The second and third products also satisfy these conditions
by Lemma~\ref{6.8}. It follows from section 3 and separate
w*-continuity of the triple product that $C$ is a w*-closed subtriple
of $A$ and that $\psi^*$ restricted to $C$ is a w*-continuous
homomorphism onto $B$. Let $C=I \oplus K$ where $K$ denotes the
kernel. Suppose $u$ is a tripotent in $B$. Let $P$ and $P^{\perp}$ be
the projections of $C$ onto $I$ and $K$. $P(u_{\psi})$ and
$P^{\perp}(u_{\psi})$ are orthogonal tripotents that sum to $u_{\psi}$
and $\psi^{\ast}(P(u_{\psi}))=u$. By Lemma~\ref{lem:6.3},
$P(u_{\psi})=u_{\psi}+q$ where $q \perp u_{\psi}$.
Hence $q=-P^{\perp}(u_{\psi})$ which forces $q^{3}=0$.
Thus $K=0$ and $\psi^*$ is a w*-continuous isomorphism from $C$ onto
$B$.  \qed

\medskip

An immediate consequence of the proof is the following corollary.
\begin{corollary}
Retain the notation of the theorem. Then $C=\{x_{\psi}|x \in B\}$.
\end{corollary}

The next two corollaries constitute a proof of Theorem 1.

\begin{corollary} 
Suppose that $A$, $B$, $C$ and $\psi$ are as in Theorem~\ref{thm:3}.  Let $\phi$ denote the inverse of $\psi^{\ast}|C$ and let $P:A_*\rightarrow A_*$ be the linear map with $P^*=\phi\circ \psi^*$ (which exists
by the automatic weak*  
continuity of $JBW^*$-triple isomorphisms).
Then $P$ is a contractive
projection of $A_{\ast}$ onto $\psi(B_{\ast})$
\end{corollary}
\pf\ For $f \in B_{*}$ and $a \in A$, $\pair{P(\psi(f))}{a}  =   
\pair{f}{\psi^{*}(\phi(\psi^{*}(a))}  =   
\pair{f}{\psi^{*}(a)}=\pair{\psi(f)}{a}$.
The statement follows.  \qed

\medskip

In the next corollary we use the following fact from the structure theory of $JBW*$-triples: every $JBW^*$-triple $U$ can be decomposed into an $\ell^\infty$-direct sum of orthogonal weak*-closed ideals $U_1$ and $U_2$, where $U_1$ is a direct sum of spaces of the form $L^\infty(\Omega,C)$, with $C$ a Cartan factor, and $U_2$ has no abelian tripotents (see \cite[(1.16)]{HorNeh88} and \cite[(1.7)]{HornMZ}).   In particular, since Hilbert spaces are Cartan factors, we can write $B=B_1\oplus B_2$ where $(B_{1})_{\ast}$ is an $\ell^1$ direct sum of
spaces isomorphic to $L^{1}(\Omega_{\lambda},H_{\lambda})$, where $H_\lambda$ is a Hilbert space of dimension at least two, and
$(B_{2})_{\ast}$ has no nontrivial $\ell^1$-summand of the from
$L^{1}(\Omega,H)$, with $H$ is a Hilbert space of dimension at least two.

\begin{corollary}
Suppose that $A$ and $B$ are JBW*-triples and $\psi$
is an isometry from $B_{\ast}$ into $A_{\ast}$, and let $B=B_1\oplus B_2$ be the decomposition described above.   Then there is a contractive projection $P$ from $A_{\ast}$ onto
$\psi((B_{2})_{\ast})$ which annihilates $\psi((B_{1})_{\ast})$
\end{corollary}
\pf\ Denote by $\psi_{i}$ the restriction of $\psi$ to
$(B_{i})_{\ast}$. It is immediate from the previous corollary that
there exists a contractive projection $P$ from $A_{\ast}$ onto
$\psi_{2}((B_{2})_{\ast})$ with $P^*=\phi_{2} \circ
\psi_{2}^{\ast}$ .  Suppose $f \in \psi_{1}((B_{1})_{\ast})$.
Pick a tripotent $u \in B_{2}$. Using Lemmas~\ref{lem:4.2} and
~\ref{lem:6.1},
$$u_{\psi_{2}}=\phi_{2}(u)=\phi_2(\sup_\lambda
v_{g_\lambda})=\sup_\lambda\phi_2(v_{g_\lambda})=\sup_\lambda
v_{\psi_2(g_\lambda)}$$
for a family of  pairwise orthogonal normal functionals $g_\lambda\in (B_2)_*$ (see
the proof of Lemma~\ref{lem:4.1}).
Since $f \perp \psi_{2}(g_{\lambda})$,  $f(v_{\psi_{2}(g_{\lambda})})=0$
and so by \cite[(3.23)]{HornMS}
$f(u_{\psi_2})=0$. Hence $f(\phi_{2}(u))=0$. It follows
that $f(\phi_{2}((\psi_{2})^{*}(A)))=0$ and $P(f)=0$.
\qed

\bibliographystyle{amsplain}

\end{document}